\documentclass[11pt,a4paper,reqno]{amsart}
\usepackage{amsfonts,amsthm,amscd,epsfig,psfrag,amsmath,amssymb,enumerate}
\numberwithin{equation}{section}
\newtheorem{theorem}{Theorem}[section]
\newtheorem{lemma}[theorem]{Lemma}
\newtheorem{proposition}[theorem]{Proposition}
\newtheorem{corollary}[theorem]{Corollary}

\newtheorem{definition}[theorem]{Definition}

\newtheorem{claim}{Claim}

\newcommand{\cB}{\ensuremath{\mathcal B}}

\newcommand{\cD}{\ensuremath{\mathcal D}}
\newcommand{\cE}{\ensuremath{\mathcal E}}
\newcommand{\cF}{\ensuremath{\mathcal F}}

\newcommand{\cL}{\ensuremath{\mathcal L}}

\newcommand{\bbE}{{\ensuremath{\mathbb E}} }

\newcommand{\bbN}{{\ensuremath{\mathbb N}} }

\newcommand{\bbP}{{\ensuremath{\mathbb P}} }

\newcommand{\bbR}{{\ensuremath{\mathbb R}} }

\newcommand{\bbT}{{\ensuremath{\mathbb T}} }

\newcommand{\bbZ}{{\ensuremath{\mathbb Z}} }

\newcommand{\eps}{\epsilon} \newcommand{\var}{{\rm Var}}
\newcommand{\be}{\begin{equation}}
\newcommand{\la}{\label} 
 \newcommand{\si}{\sigma}

\newcommand{\wt}{\widetilde}
\newcommand{\wb}{\bar}

\newcommand{\un}{{\ensuremath{\bf 1}}}

\DeclareMathSymbol{\leqslant}{\mathalpha}{AMSa}{"36} % nicer `smaller or equal'
\DeclareMathSymbol{\geqslant}{\mathalpha}{AMSa}{"3E} % nicer `larger or equal'
\DeclareMathSymbol{\eset}{\mathalpha}{AMSb}{"3F}     % nicer `emptyset'
\renewcommand{\leq}{\;\leqslant\;}                   % redef. of < or =
\renewcommand{\geq}{\;\geqslant\;}                   % redef. of > or =
\newcommand{\cgap}[0]{c_{\mathrm{gap}}}
\newcommand{\csob}[0]{c_{\mathrm{sob}}}
\newcommand{\df}[0]{{\mathcal{D}}}

\newcommand{\Cov}[0]{{\mathrm{Cov}}}

\def\Tree{{\mathbb T}}
\def\muAeta{{\mu_A^\eta}}
\def\VarAeta{{\var_A^\eta}}
\def\EntAeta{{\Ent_A^\eta}}

\let\neper=e

\def\ie{\hbox{\it i.e.\ }}

\let\sset=\subset

\def\nep#1{ \neper^{#1}}

\def\tc{\thinspace | \thinspace}
\let\<=\langle
\let\>=\rangle

\def\Ent{ \mathop{\rm Ent}\nolimits }

\def\ninf#1{ \| #1 \|_\infty }

\def\inte#1{\lfloor #1 \rfloor}

\def\smallno{\smallskip\noindent}
\def\medno{\medskip\noindent}

\def\\{\hfill\break}

\def\Todd{\Tree_{\rm odd}}
\def\Teven{\Tree_{\rm even}}

\let\a=\alpha \let\b=\beta   \let\d=\delta  \let\e=\varepsilon
 \let\g=\gamma \let\h=\eta      \let\l=\lambda
\let\m=\mu      \let\o=\omega      
\let\r=\rho  \let\s=\sigma \let\t=\tau   
  
\let\D=\Delta   \let\G=\Gamma  \let\L=\Lambda 
\let\O=\Omega      

\begin{document}

%% \title[Glauber dynamics on trees] {On the stationary measures for the
%% stochastic Ising model and hard core gas on a regular tree}
\title[Phase ordering after a deep quench]{Phase ordering after a deep quench: the stochastic Ising and
  hard core gas models on a tree}
\date{\today} 
\author[P. Caputo]{Pietro Caputo}
\address{Dip. Matematica, Universita' di Roma Tre, L.go S. Murialdo 1,
00146 Roma, Italy} \email{caputo\@@mat.uniroma3.it}\thanks{}
\author[F. Martinelli]{Fabio Martinelli} \address{Dip. Matematica,
Universita' di Roma Tre, L.go S. Murialdo 1, 00146 Roma, Italy}
\email{martin\@@mat.uniroma3.it}\thanks{}

\vskip 1cm
\begin{abstract}
\noindent 
Consider a low temperature stochastic Ising model in the phase
coexistence regime with Markov semigroup $P_t$. A fundamental and still
largely open problem is the understanding of the long time behavior of
$\d_\h P_t$ when the initial configuration $\h$ is sampled from a highly
disordered state $\nu$ (e.g. a product Bernoulli measure or a high
temperature Gibbs measure).  Exploiting recent progresses in the
analysis of the mixing time of Monte Carlo Markov chains for discrete
spin models on a regular $b$-ary tree $\Tree^b$, we tackle the above
problem for the Ising and hard core gas (independent sets) models on
$\Tree^b$. If $\nu$ is a biased product Bernoulli law then, under
various assumptions on the bias and on the thermodynamic parameters, we
prove $\nu$-almost sure weak convergence of $\d_\h P_t$ to an extremal Gibbs
measure (pure phase) and show that the limit is approached at
least as fast as a stretched exponential of the time $t$. In the
context of randomized algorithms and if one considers the Glauber
dynamics on a large, finite tree, our results prove fast local
relaxation to equilibrium on time scales much smaller than the true
mixing time, provided that the starting point of the chain is not taken
as the worst one but it is rather sampled from a suitable distribution.
\end{abstract}

\maketitle

\section{Introduction}
Let $G=(V,E)$ be a countable infinite graph of
bounded degree and consider, for definiteness, 
a continuous time stochastic
Ising model (Glauber dynamics) $\bigl\{\s_t^\h\bigr\}_{t\geq 0}$ on $G$ with
initial condition $\h$ and infinitesimal generator $\cL$. Here $\eta$
is picked from the set $\O$ of assignments of a \ $\pm 1$
variable to each vertex $x\in V$. The main problems discussed in this 
paper can be formulated as follows.

Assume that the thermodynamic parameters 
are such that there exist multiple reversible Gibbs
measures for $\cL$. For stochastic
Ising models this amounts to say that the inverse temperature $\b$ and
the external field $h$ are such that $\mu^+ \neq \mu^-$, where $\mu^+$
and $\mu^-$ are the Gibbs measures obtained by taking infinite
volume limits with pure $+$ and $-$ boundary conditions, respectively.
Suppose that $\h\in\O$ is distributed according to a Bernoulli product
measure with parameter $p$, i.e.\ $\{\eta_x\}_{x\in V}$ is a
collection of i.i.d.\ random variables with $\bbP(\h_x=+1)=p$. 
%In other
%words $\h$ is chosen according to an infinite temperature Gibbs measure
%with density $p$ of plus spins. The relevant questions are then the
%following. 
Then:
\begin{enumerate}[i)]
\item  Under which condition on the bias $p$ is the Ising
plus phase $\mu^+$ the unique limit point of the law of 
$\s_t^\h$ as $t\to\infty$, for
a.a.\ $\h$ ? 
\item If so, how fast does the law of $\s_t^\h$
approach $\mu^+$ ?
\end{enumerate}

The above questions, with $G$ some regular lattice, 
have their origin in the theory of ``phase ordering
kinetics'' \cite{Bray}  -- that is growth of order through a dynamical domain coarsening 
--  and clearly represent basic problems in the theory of
interacting particle systems. Unfortunately, a rigorous approach these problems is still
largely missing.
%
%Unfortunately the problem of the characterization of the invariant
%measures for stochastic Ising models is still largely open. \\ 
%When $G=\bbZ^d$ it is known that every invariant measure for $\cL$ which
%is also translation invariant is a Gibbs measure and in the special case
%$d=2$ the same conclusion hold without the translation invariant
%hypothesis \cite{Li}. 

If the law of the starting configuration $\eta$ stochastically
dominates the plus phase $\mu^+$ it is possible to use some
monotonicity arguments 
(allowed by the ferromagnetic character of the model)
to prove that $\mu^+$ is indeed the unique
limiting point of the process and that the convergence takes place
faster than any inverse
power of $t$. We refer to section 6.7 in \cite{Mar} for the case
$G=\bbZ^d$ and to our Lemma \ref{mu+le} below for a
stronger statement in the case of regular trees. 
However, it is easily seen that such a stochastic domination
requirement for the initial Bernoulli distribution forces the bias $p$
to be exponentially close to $1$ when $\b\to\infty$. For $\b<\infty$ we do not
know of any result that goes beyond this
simple case.
 
On the other hand, the extreme case $\b=\infty$ (zero temperature
Glauber dynamics) has received considerable attention in the
probabilistic literature, and various kinds of graphs ($\bbZ^d$, the hexagonal
lattice and the binary tree) have been considered
\cite{CDN,HN,CNS,H,FoSchSi}. In this case,
besides the motivation from physics to study simple models of spatial
domain coarsening, there is also an interesting connection with
(non--linear) voter models \cite{Li2}.  The relevant quantities are then
the probability that a given vertex flips its value finitely or
infinitely many times, the probability that a given spin has not flipped
before time $t$, the typical size of clusters of vertices with a common
spin value and other related percolation
questions.

Going back to our original problems, 
a major obstacle for
progresses in the case $G=\bbZ^d$ is represented by the absence of tight
bounds on the mixing time of the Glauber dynamics in finite boxes with
plus boundary conditions, \ie those boundary conditions that select
the plus phase. On the contrary, when $G$ is the regular $b$--ary
tree,
this question has been recently solved in a sharp and constructive
way for various models \cite{MaSiWe}. Exploiting the results of 
\cite{MaSiWe} we have been able to study
the above basic questions  
for two attractive systems on trees: the Ising model and the hard core
gas (independent sets). Our results provide some answers to i) and ii)
in non trivial cases. For instance we show that if the bias is
sufficiently large (but independent of $\b$) then we have the desired
convergence for {\em all} temperatures. The paper also includes
a discussion of several interesting problems
that are left unsolved and that we would like to consider in future work. 
For simplicity we present now our main result only for the Ising model,
and defer the reader to section 6 for the analogous theorem for the
hard-core gas. Before stating our results we will now briefly overview
the model and its basic features.

\subsection{The Ising model on the $b$-ary tree}
\label{sec:Ising_on_trees}
\begin{figure}[h]
\centerline{\psfig{file=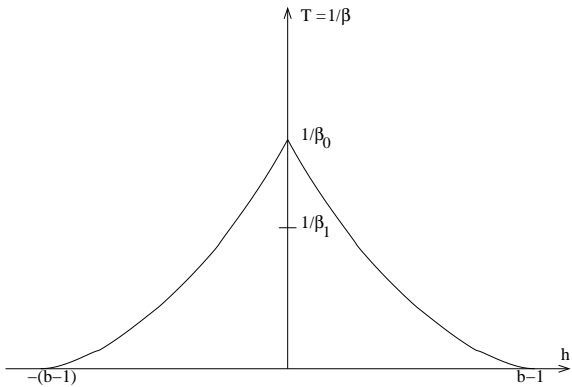,height=2in}}
\caption{The critical field $h_c(\beta,b)$.  The Gibbs
measure is unique above the curve.}
\label{fig:phase_diagram}
\end{figure}

%\vskip-0.05in 
From now on $\bbT^b$ denotes the infinite, rooted $b$--ary tree,
where each vertex has exactly $b$ children ($b\geq 2$ is a given integer).
The Ising Gibbs measure on~$\Tree^b$  at inverse
temperature $\b$ and external field $h$, formally given by 
\begin{equation*}
 \mu(\s)\propto \exp\left[\b\bigr(\sum_{xy\in E}\s_x\s_y + h\sum_x
 \s_x\bigr)\right]\,, %\quad \s\in \O:=\{-1,1\}^{\Tree^b} 
\end{equation*}
where $E$ is the set of edges of $\Tree^b$, has recently received a lot
of attention as the canonical example of a statistical physics model on
a ``non-amenable'' graph (i.e., one whose boundary is of comparable size
to its volume) -- see e.g.
\cite{BRZ,Ioffe,EKPS,ST98,JS99,BKMP,BRSSZ}. 
The phase diagram of the model in the $(h,\b)$ plane is
known (\cite{Georgii,L}) to be 
quite different from that on the cubic lattice $\bbZ^d$ (see
Fig.~\ref{fig:phase_diagram}).

We now recall some of its basic features. We write $T_\ell$ for the 
rooted tree obtained by removing all vertices which are at distance
greater than $\ell$ from the root. The measures $\mu^+$ and $\mu^-$
are obtained by imposing $+1$ and, respectively, $-1$ boundary data at
the leaves of $T_\ell$ and taking the limit $\ell\to\infty$. The free
measure $\mu^{\rm free}$ is defined as the limit $\ell\to\infty$ when
the boundary data at the leaves of $T_\ell$ are free (\ie absent).

On the line $h=0$ there is a first critical
value \hbox{$\b_0= \frac12 \log\bigl(\frac{b+1}{b-1}\bigr)$}, marking the dividing
line between uniqueness and non-uniqueness of the Gibbs measure (\ie
$\mu^+\neq\mu^-$
as soon as $\b\geq \b_0$). Then, in sharp
contrast to the model on~$\bbZ^d$, there is a second critical point
\hbox{$\b_1= \frac 12\log\bigl(\frac{\sqrt{b}+1}{\sqrt{b}-1}\bigr)$} which is often
referred to as the ``spin-glass critical point''~\cite{CCST} and has
different interpretations.
If one considers for instance the model with $h=0$ on the finite tree
$T_\ell$ with
i.i.d.\ Bernoulli random boundary data $\eta$ with $p=1/2$ at the leaves of $T_\ell$,
then the distribution of the magnetization at the root (as a function of
$\eta$) becomes non
trivial only if $\b>\b_1$, see \cite{CCST}. In particular, as $\ell \to
\infty$, for $\b\leq \b_1$ the Gibbs measure on $T_\ell$ with the
above random boundary $\eta$,
converges (weakly) a.s.\ to the free measure $\mu^{\rm free}$. 
Another way to look at $\b_1$ is to say that $\mu^{\rm free}$ is an extremal
Gibbs measure iff $\b\leq\b_1$ (see ~\cite{BRZ,Ioffe,Ioffe2,BKMP} and,
more recently, \cite{MaSiWe}). Finally $\b_1$ has also the interpretation of
the non-reconstruction/reconstruction threshold in the context of ``bit
reconstruction problems'' on a noisy symmetric channel \cite{EKPS,Mossel,MP01}.

When an external field $h$ is added to
the system, it turns out that for all $\beta>\beta_0$, there is a
critical value~$h_c=h_c(\beta,b)>0$ of the field such that
$\mu^+\neq\mu^-$ iff $|h|\leq h_c$. 
The Ising model on the tree at external field $h=\pm h_c$ therefore
shares the following two properties with the classical 
Ising model on $\bbZ^d$ at zero external field: on one hand the Gibbs
measure is sensitive to the choice of boundary condition; on the other hand
any arbitrarily small increase of $|h|$ causes the Gibbs measure to
become insensitive to the boundary condition. 

\subsection{The Glauber dynamics}
\label{sec:dynamics}
%\vskip-0.05in
The {\it Glauber dynamics\/} on $\Tree^b$ is the unique Markov process
$\{\s_t^\h\}_{t\geq 0}$
on~$\Omega$ with $\s_{t=0}^\h=\h$ and Markov generator $\cL$ formally given by
\begin{equation}
  \label{generator}
  (\cL f)(\s) = \sum_{x\in \bbT^b}c_x(\s) [f(\s^{x})-f(\s)]\,,
\end{equation}
where $\s^x$ denotes the configuration obtained from $\s$ by flipping
the spin at $x$, and $c_x(\s)$ denotes the flip rate at $x$.

Glauber dynamics on trees has received recently considerable interest
~\cite{BKMP,MaSiWe}. 
Results in \cite{MaSiWe} show in a
rather strong form that the mixing time (see
e.g.\ \cite{Saloff} for a definition) 
on the finite subtree $T_\ell$ is always
$O(\ell)$ if either $\b<\b_1$ and $h$ is arbitrary or if $\b,h$ are
arbitrary and the boundary conditions on the leaves of $T_\ell$ are identically equal
to $+1$ (or, by symmetry, to $-1$). 
In particular, the Glauber dynamics in the pure plus phase $\mu^+$
always mixes fast (see \cite{MaSiWe} and section 2 below for more details).   

Although all our results apply to any choice of finite--range, 
uniformly positive, bounded and attractive  flip rates satisfying
the detailed balance condition w.r.t.\ the Ising Gibbs
measure (see \cite{Mar}), for simplicity in the sequel we will work with a specific 
choice known as the {\it heat-bath\/} dynamics (see section 2 below
for the definition).
We will use the standard notation $P_t=\nep{t\cL}$ for
the Markov semigroup associated to $\cL$. The spin at $x$ at time $t$
with starting configuration $\eta$ is denoted by $\si_{t,x}^\eta$ and
we will often use the shortcut notation
\be
\r_{t,x}(\eta) %=\bbE(\si_{t,x}^\eta) 
= (P_t\si_x)(\eta)
\la{rhodef}
\end{equation}
for the expected value $\bbE(\si_{t,x}^\eta)$ 
of $\si_{t,x}^\eta$ given that the process starts in $\eta$.

\subsection{Main results}
In order to state our main results we need an extra bit of notation. 
We first define the set of initial configurations $\eta$ such that the Glauber
dynamics $\s_t^\eta$ converges weakly, at a certain rate, to the plus phase $\mu^+$. 
\begin{definition}
\label{Omega alpha}
Given $\a\in (0,1)$, $\O_\a$ will denote the set of
starting configurations $\eta\in\O$ such that for any $x\in \Tree^b$ there exists a
time $t_0=t_0(\eta,x)<\infty$ such that for all $t\geq t_0$
\begin{equation}
|\,\rho_{t,x}(\eta)- \mu^+(\s_x)\,|\leq \exp(-t^\a)\,.
\label{1}
\end{equation}
\end{definition}
We will see in Corollary \ref{weak} below 
that for any $\eta\in \O_\a$ the law of the process $\s^\eta_t$ converges weakly to
$\mu^+$ as $t\to \infty$.

The initial configuration $\eta$ is often sampled from a 
Bernoulli product measure with parameter $p$, \ie $\eta_x
= +1$ with probability $p$ and $\eta_x = -1$ with probability $1-p$
independently for each $x\in \Tree^b$. We write $\bbP_p, \bbE_p$ for the
corresponding probability and expectation. 
Finally, we need to recall the notion of partial ordering 
(stochastic domination) between
probability measures on $\O$. Given two configurations $\s,\h\in \O$ we
will write $\s\leq\h$ iff $\s(x)\leq\h(x)\; \forall \, x\in \Tree^b$. A
function $f:\O\mapsto\bbR$ is called {\it monotone increasing
  (decreasing)\/} if $\s\leq\s'$ implies $f(\s)\leq f(\s')$ ($f(\s)\ge
f(\s')$). Given two probability measures $\mu$, $\mu'$ on $\O$ we write
$\mu \leq \mu'$ if $\mu(f) \leq\mu'(f)$ for all (bounded and measurable)
increasing functions $f$.
\smallno

Our main results can now be stated as follows. 
\begin{theorem}
\label{theotree}
\medno
\begin{enumerate}[a)]
\item For every $a>0$, $b\geq 2$,
there exists $p<1$ such that
for all $\b\in(0,\infty)$ and $h\geq -h_c(\b,b) + a$ 
we have $\nu(\O_\a)=1$, for some $\a=\a(\b,h,b)>0$, for
any initial distribution 
$\nu$ such that $\nu\ge
\bbP_{p}$\,.
\medno
\item For every $p>\frac12$, there exist $b_0\in\bbN$ and $\b_0\in(0,\infty)$ such that 
for $h=0$, $b\geq b_0$, $\b\geq \b_0$ we have $\nu(\O_\a)=1$, for some $\a=\a(\b,b,p)>0$,
for any initial distribution $\nu$ such that $\nu\geq\bbP_{p}$\,.
\medno
\item Let $\O^-_\a$ denote the event defined in (\ref{1}) with $\mu^-$
in place of $\mu^+$. For every $p<1$ there exist
$b_0\in\bbN$ and $\b_0\in(0,\infty)$ such that 
for $b\geq b_0$, $\b\geq \b_0$ and $h=-h_c(\b,b)$,
we have $\nu(\O^-_\a)=1$ for some $\a=\a(\b,b,p)>0$, for all 
initial distributions $\nu$ such that $\nu\leq \bbP_p$.
\end{enumerate}
\end{theorem}

\smallskip

\subsection{Remarks}
Let us make some remarks on the above statements. 

\medno
{\bf 1}. We believe that in the case $h=0$,
convergence to the plus phase should occur as soon as $p>\frac12$. 
Unfortunately our bounds on $p$ in statement \emph{a)} are far from being sharp. 
However, as stated in \emph{b)}, we can approach the critical 
value $\frac 12$, by taking $b$ large. 
Another interesting issue is the dependence of $p$ on $h$. Our technique in 
the proof of part a) of Theorem \ref{theotree} breaks down in the case $h=-h_c(\b,b)$
and the value of $p$ in that statement approaches 
$1$ as $a\to 0$. On the other hand, statement \emph{c)} shows that if $h=-h_c(\b,b)$
the critical value of $p$ for convergence to the plus phase must approach $1$
when $b\to\infty$.

\medno
{\bf 2}. The main arguments we use to prove Theorem \ref{theotree}
are based on two essential features
of the Ising model on $\bbT^b$. 
The first is monotonicity which is shared
by all so--called {\em attractive} interacting particle systems.
The second is the so--called ``rigidity'' of critical phases for
spin systems on trees (\cite{BRSSZ}). Roughly speaking the latter
means that, as long as $h > - h_c$, if we are in the pure phase 
$\mu^+$ we can add a small density of spins of the
opposite ($-$) phase and this will not alter significantly the structure
of $\mu^+$. This, in a sense, is what we do when we introduce
{\em obstacles} (see section 3 below) to lower
bound the magnetization $\r_{t,x}(\eta)$. 
The hard core gas model will be shown to have
both these properties and our results there
(see Theorem \ref{hard-corethm} below) will be obtained essentially by the same
methods. On the other hand these techniques do not apply when there is
no rigidity of phases, as e.g.\ in the Ising model on $\bbZ^d$
(see \cite{ScSh} for a deep investigation of the metastable behavior
of this model), or
when there is no attractivity, as e.g.\ in the $q$--state Potts model for
$q\geq 3$.  

\medno
{\bf 3}. A close check of the various probabilistic estimates needed for the
proof of Theorem \ref{theotree} and which are behind a Borel--Cantelli
characterization of the set $\O_\a$, shows that there is also an
$L^2$-version of Theorem \ref{theotree}, with (\ref{1}) replaced
by a bound of the form:
\begin{equation*}
\bbE_p\bigl((\rho_{t,x}^\h- \mu^+(\s_x))^2\bigr)\leq \exp(-t^\a) \,,
\end{equation*}
for any $t$ large enough.

\subsection{Plan of the paper}
The rest of the paper is organized as follows. In section 2 we
give the basic preliminaries for the proof of Theorem \ref{theotree}.
In section 3 we explain our main argument. In particular, 
here we give the proof of Theorem \ref{theotree} by
assuming the validity of several technical claims.
Section 4 and 5 deal with the proof of these claims. 
In section 6 we present our results for the hard core gas. 
Finally, some further results and open problems will be discussed
in section 7. 

\section{Some preliminaries}
\label{sec:prelims}
%\vskip-0.3in\hbox{}
Here we first collect several useful preliminaries
concerning the Gibbs measure and the Glauber dynamics and then 
discuss some basic results on convergence to the plus phase,
together with properties of the sets $\O_\a$ introduced
above.

\subsection{Finite Gibbs measures on the $b$-ary tree}
%\vskip-0.05in
%For $b\geq 2$, let $\Tree^b$ denote the infinite, rooted $b$-ary tree (in
%which every vertex has $b$ children). 
We denote by $d(x,y)$ the tree
distance between two vertices $x,y\in\bbT^b$. If $r$ is the root of the
tree, we write $d(x)=d(x,r)$ for the depth of $x$.  When $A$ is a subset
of vertices of $\bbT^b$ we set $d(x,A)=\inf_{z\in A}d(x,z)$. The
boundary of $A$, $\partial A$, is defined as the set of vertices $x$
such that $d(x,A)=1$. $E(A)$ denotes the set of $\bbT^b$--edges $(x,y)$
with $x,y\in A$. 

The Ising spin configurations space is the set $\O=\{-1,+1\}^{\Tree^b}$
and its elements will be denoted by Greek letters $\s,\eta,\xi$ etc. The
set $\O$ is equipped with the standard $\s$--algebra $\cF$
generated by the variables $\{\s_x\}_{x\in \Tree^b}$.  For any finite
subset $A\subseteq \Tree^b$ and any $\h\in \O$, we denote
by~$\mu_A^\eta$ the Gibbs distribution over~$\Omega$ conditioned on the
configuration outside~$A$ being~$\eta$: i.e., if $\sigma\in\Omega$
agrees with~$\eta$ outside~$A$ then
$$
\mu_A^\eta(\sigma) \propto
\exp\Bigl[\beta\bigl(\sum\nolimits_{xy\in E(A\cup\partial
  A)}\sigma_x\sigma_y + h\sum\nolimits_{x\in A}\sigma_x\bigr)\Bigr], 
$$
where $\beta$ is the inverse temperature and $h$ the external field. We
define $\mu_A^\eta(\sigma)=0$ otherwise. If the boundary configuration
$\eta$ is identically equal to $+1$ ($-1$) we will denote the
corresponding conditional Gibbs distribution by $\mu_A^+$
($\mu_A^-$). Whenever the set $A$ will coincide with the finite subtree $T_\ell=\{x\in
\Tree^b:\;d(x)\leq \ell\}$ we will abbreviate the symbol $T_\ell$ in the Gibbs
measure with $\ell$, \ie $\mu_{\ell}^\eta$ stands for $\mu_{T_\ell}^\eta$.

For a bounded measurable function $f:\Omega\to\bbR$ we denote by
$\mu_A^\eta(f)= \sum_{\sigma\in\Omega} \mu_A^\eta(\sigma)f(\sigma)$ the
{\it expectation\/} of~$f$ w.r.t.\ the distribution~$\mu_A^\eta$.
Analogously, for any $X\in \cF$, $\mu_A^\eta(X):= \mu_A^\eta(\un_X)$ where
$\un_X$ is the characteristic function of the event $X$.  We will write
$\var_{\muAeta}(f)$ or $\VarAeta(f)$ for the \emph{variance}
$\muAeta(f^2)-\muAeta(f)^2$ and (for $f\geq 0$) $\Ent_{\muAeta}(f)$ or
$\EntAeta(f)$ for the \emph{entropy} $\muAeta(f\log
f)-\muAeta(f)\log\muAeta(f)$ w.r.t.\ $\muAeta$.  Note that
$\VarAeta(f)=0$ iff, conditioned on the configuration outside~$A$
being~$\eta$, $f$~does not depend on the configuration inside~$A$.  The
same holds for $\EntAeta(f)$. We shall use the 
symbol $\mu_A$ for the map $\eta\to\mu_A^\eta$. Similarly,
$\var_{A}$ and $\Ent_A$ stand for $\eta\to\var_{A}^\eta$ and
$\eta\to\Ent_A^\eta$. 

%\begin{definition}
A probability measure $\mu$ on $\bigl(\O,\cF\bigr)$ will be called a
Gibbs measure for the Ising model with parameters $(\b,h)$ if 
$$
\mu\bigl(\mu_A(X)\bigr)=\mu(X),\qquad \text{ for all }X\in \cF \ \text{
  and all finite sets } A\sset \Tree^b\,.
$$
%\end{definition}
In this work a crucial role will be played by the following monotonicity
property of the Gibbs measures (and of the Glauber dynamics, see below) known
as \emph{attractivity}. For any
increasing bounded measurable function $f$:
\begin{align}
&(i) \qquad \text{for any $A\sset \Tree^b$ the map $\h \mapsto
  \mu_A^\eta(f)$ is increasing;}
\label{mono1}\\
&(ii) \qquad \text{$\mu_B^+(f)\leq \mu_A^+(f)$ whenever $A\sset B$.}
  \label{eq:mono20}
\end{align}
Recall that 
%Among the Gibbs measures (if more than one) for the Ising model in this
%paper we shall mostly consider 
the ``plus phase'' $\mu^+$ is obtained as
the weak limit as $\ell\to\infty$ of $\mu^+_\ell$. 
Existence of this limit follows from the
monotonicity properties described above. Similarly one defines the
``minus phase'' $\mu^-$. It turns out that any 
(infinite volume) Gibbs measure $\mu$ satisfies
$\mu^-\leq \mu\leq \mu^+$.

\subsection{The Heat Bath dynamics on finite trees}
\label{sec:glauber}
%\vskip-0.05in
For any finite
subset $A\subseteq \Tree^b$ and any $\t\in \O$
we define the Heat Bath Glauber
dynamics in $A$ with boundary condition (b.c.)
$\t$ (see e.g \cite{Mar}) as the continuous time Markov chain on
$\O_A^\t:=\{\s\in \{-1,1\}^{A\cup\partial A}:\ \si=\t \;{\rm on}\;{\partial A}\}$
with generator
\begin{equation}
  \label{generatorA}
  (\cL_A^\t f)(\s) = \sum_{x\in A}c_x(\s) [f(\s^{x})-f(\s)]\,, \quad
  \s\in \O_A^\t\,,
\end{equation}
where $(\s^x)_y=\si_y$ for all $y\neq x$ and $(\si^x)_x=-\si_x$ 
and
$$
c_x(\s)= \mu^{\s}_{\{x\}}(\s^x) = \frac{1}{1+w_x(\sigma)}, \quad 
%\text{where}\quad
w_x(\sigma):=\exp\left[2\beta\sigma_x\left(h+\sum\nolimits_{y:\ d(x,y)=1}\sigma_y \right)\right].
$$
In analogy with the infinite volume case discussed
in the introduction the chain started from $\xi$
will be denoted by $\{\s_t^{\xi,A,\t}\}_{t\geq 0}$. If $A=T_\ell$ 
we will simply write $\s_t^{\xi,\ell,\t}$. 

It is well known that there is a global
pathwise coupling among the processes $\bigl\{(\s_t^{\xi;A,\t})_{t\ge
  0},\ A\sset \Tree^b,\,\xi,\t\in\O\bigr\}$ such that, for any $A\sset
B\sset \Tree^b$, any $\xi\leq \xi'$ and any $\t\leq \t'$:
\begin{align}
  \s_t^{\xi;A,\t} &\leq \s_t^{\xi';A,\t'}  &\forall t\geq 0 \nonumber\\
\s_t^{\xi;A,-} &\leq \s_t^{\xi;B,\t} \leq \s_t^{\xi;A,+}  &\forall t\geq 0
\label{eq:mono2}
\end{align}
It is a well--known (and easily checked) fact that, for any finite
$A\sset \Tree^b$ and any $\t$, the Glauber 
dynamics in $A$ with b.c.\ $\t$ is ergodic and reversible w.r.t.\ the Gibbs 
distribution~$\mu_A^\t$, \ie for any function $f$  
$$
\lim_{t\to \infty}\nep{t\cL_A^\t}f =\mu_A^\t(f)\,.
$$
The rate at which the above convergence takes place 
is often measured using two concepts from
functional analysis: the {\it spectral gap\/} and the {\it logarithmic
  Sobolev constant}. We now describe these two quantities for a generic
(finite or infinite volume) Gibbs measure $\mu$.

For a local function $f:\Omega\to\bbR$ define the
{\it Dirichlet form\/} of~$f$ associated to the Glauber dynamics with
reversible measure $\mu$ by
\begin{equation}
\label{eq:df}
\df_\mu(f)= \frac 12 \sum_x \mu\bigl(c_x \bigl[f(\s^x)-f(\s)\bigr]^2\bigr)
   = \sum_x\mu(\var_{\{x\}}(f)).
\end{equation}
(The l.h.s.\ here is the general definition for any choice of the
flip rates~$c_x$;
the last equality holds when specializing to the case of the heat-bath
dynamics.) 
The \emph{spectral gap} $\cgap(\mu)$ and the \emph{logarithmic Sobolev
constant} $\csob(\mu)$ of $\mu$ are then
defined by 
\begin{eqnarray}
   \cgap(\mu) = \inf_{f} {{\df_\mu(f)}\over\var_\mu(f)};\qquad
   \csob(\mu) = \inf_{f\geq 0} {{\df_\mu(\sqrt{f}\,)}\over\Ent_\mu(f)},\label{eq:cgapcsob}
\end{eqnarray}
where the infimum in each case is over non-constant functions~$f$.

The spectral gap $\cgap(\mu)$ measures the rate of exponential decay
as $t\to \infty$ of the variance w.r.t.\ $\mu$, \ie $ \cgap(\mu)$ is
the (largest) constant such that for any $f$
\be
\var_\mu(P_tf)\leq \nep{-2t\cgap(\mu)}\,\var_\mu(f)\,,
\la{gapvar}
\end{equation}
where $P_t$ denotes the semigroup associated to the Dirichlet form $\df_\mu(f)$.
The log--Sobolev constant $\csob(\mu)$ is related to the following
hypercontractivity
estimate (see e.g.\ ~\cite{Saloff}): setting
$q_t:=1+\nep{4\csob(\mu)t}$ we have, for any function $f$ and any
$t\geq 0$ 
\be
\|P_t f\|_{q_t,\mu}\leq \|f\|_{2,\mu}\,,
\la{hypp}
\end{equation} 
where $\|f\|_{p,\mu}$ stands for the $L^p$--norm of $f$ w.r.t.\ $\mu$.

If $\mu$ is a \emph{finite volume}
Gibbs measure (\ie $\mu=\mu_A^\t$) then both $\cgap(\mu)$ and
$\csob(\mu)$ are always strictly positive (possibly depending on
$A,\t$). The striking result of \cite{MaSiWe} is that the same is true for
\emph{any} choice of the parameters $(\b,h)$ if $\mu=\mu^+$ is the
infinite volume plus phase. More
precisely one has
\begin{equation}
  \label{eq:key-bound}
  \inf_{\ell} \cgap(\mu_\ell^+) >0,\quad \inf_{\ell} \csob(\mu_\ell^+) >0 \,.
\end{equation}
Such a
result does not imply however any ergodicity statement for the infinite
volume Glauber dynamics. Simple monotonicity considerations show in fact
that for any increasing local function $f$ and any $t\geq 0$:
$$
%\nep{t\cL}
P_tf(-)\leq \mu^-(f)\leq \mu^+(f)\leq P_tf(+) \,, 
$$
\ie non--ergodicity whenever $\mu^-\neq \mu^+$, that is $\b>\b_0$ and $|h|\le
h_c(\b)$.

\subsection{Convergence to the plus phase: preliminary results}
A first important step in the proof of Theorem \ref{theotree}
is to show that convergence to $\mu^+$ occurs when we start from all
$+$ spins. Recall that $\r_{t,x}(\eta) = \bbE(\si^\eta_{t,x})$
stands for the expectation at time $t$ under the infinite--volume dynamics
started in $\eta$.   
\begin{lemma}
\label{pro+}
For all $b,\b,h$ there exist
$\d>0$ such that the following holds.
For all $x\in \Tree^b$ there exists $t_0(x)<\infty$ such that 
if $t\geq t_0(x)$ then  
\begin{equation}
0\leq \r_{t,x}(+) - \mu^+(\s_x) \leq \exp(-\d \,t)
\label{2}
\end{equation}
\end{lemma}
\begin{proof}
The left inequality is a direct consequence of monotonicity 
(see (\ref{eq:mono2})) and the fact that 
$\mu^+(\r_{t,x})=\mu^+(\s_x)$ for all $x$ since $\mu^+$ is an
invariant measure. 
We now prove the right inequality. For
simplicity we only analyze the case of the root $x=r$ (the general case 
requires no modifications in the argument.)
Fix a length scale $\ell$ and observe that for any $\eta$ 
monotonicity implies $\r_{t,r}(\eta)\leq
\r_{t,r}^{\ell,+}(\eta)$, with the latter denoting
expectation of $\si_{t,r}^{\eta,\ell,+}$ (the spin at the root at time
$t$ for the dynamics in $T_\ell$ with $+$ b.c.\ at the leaves of
$T_\ell$ and initial condition $\eta$).  
Let also $\phi_t^{\ell,+}$ denote the function $\eta\to 
\r_{t,r}^{\ell,+}(\eta) - \mu_\ell^+(\s_r)$ so that 
$$
\r_{t,r}(+) - \mu^+(\s_r) \leq 
\left[\mu_\ell^+(\s_r)-\mu^+(\s_r)\right] \,+\, \phi_t^{\ell,+}(+) \,.
$$
%From the uniform bound on the
%logarithmic Sobolev constant (\ref{eq:key-bound}), 
Setting $q_t:=1+\nep{2\csob(\mu_\ell^+) t}$, the estimate (\ref{hypp})
yields 
\begin{equation*}
\| \phi_t^{\ell,+}\|_{q_t,\mu_{\ell}^+}
\leq \|\phi_{t/2}^{\ell,+} \|_{2,\mu_{\ell}^+}\,.
\la{30}
\end{equation*}
Let $\mu_{\ell}^+(+)$ denote the probability of having all $+$
spins in $T_\ell$ under $\mu_{\ell}^+$. Then there exists $C_1<\infty$
such that $\mu_{\ell}^+(+)\geq \nep{-C_1b^\ell}$. Moreover by 
(\ref{eq:key-bound}) we know that $q_t\geq \nep{c_2 t}$ for some
positive $c_2$ independent of $\ell$. Then  
%An application of the uniform 
%spectral gap estimate (\ref{eq:key-bound}) then yields
\begin{align}
\phi_t^{\ell,+}(+)&\leq (\mu_{\ell}^+(+))^{-\frac1{q_t}}
\| \phi_t^{\ell,+}\|_{q_t,\mu_{\ell}^+}
\nonumber \\
& \leq \exp{(C_1b^\ell\nep{-c_2 t})} \,
\|\phi_{t/2}^{\ell,+} \|_{2,\mu_{\ell}^+}\,.
\label{3}
\end{align}
Set now $\ell=c_3\, t$
with $c_3>0$ small enough. Using (\ref{gapvar}) and
(\ref{eq:key-bound})
we therefore arrive at
$$
\phi_t^{\ell,+}(+) \leq \exp{(C_1b^\ell\nep{-c_2 t})} 
\nep {-\cgap(\mu_\ell^+)t}\leq \nep{-c_4 t}
$$
for a suitable constant $c_4>0$ and $t$ sufficiently large. 
Now the claim
(\ref{2}) follows from the fact (see e.g.\ \cite{MaSiWe}) 
that in the $+$ phase
the influence of plus boundary conditions decays exponentially fast at
any temperature: there exists $c_5>0$ such that
\begin{equation*}
|\mu^+(\s_r)-\mu_\ell^+(\s_r)| \leq \nep{-c_5 \ell}\,.
\label{4}
\end{equation*}\end{proof}
The previous result allows to show that the set $\O_\a$ is
increasing, \ie its indicator function is increasing. 

\begin{corollary}
\label{incr}
For any $\a\in (0,1)$ the event $\O_\a$ is increasing.  
\end{corollary}
\begin{proof}
We need to show 
that for any pair $(\h',\h)$ with $\h'\geq \h$ and $\h\in
\O_\a$, also the first component $\h'$ belongs to
$\O_\a$. To prove the claim we observe that, for any $x\in \Tree^b$ and any $t\ge
0$, monotonicity implies
  \begin{equation}
    \rho_{t,x}(\h)-\mu^+(\s_x) \leq \rho_{t,x}(\h')-\mu^+(\s_x)
\leq \rho_{t,x}(+)-\mu^+(\s_x)
  \label{mon}
\end{equation}
The l.h.s. of (\ref{mon}) is $\geq -\nep{-t^\a}$ for any large
enough time $t$ because $\h\in \O_\a$ by assumption.
The r.h.s. is instead bounded via Lemma \ref{pro+} above.
\end{proof}

Another consequence of Lemma \ref{pro+} is the following
\begin{corollary}
\label{weak}
For any
  $\h\in \O_\a$ the law of the process $\s^\h_t$ converges weakly to
  $\mu^+$ as $t\to \infty$.
\end{corollary}
\begin{proof} 
Observe first that by the global coupling, for any $x\in\bbT^b$ and
$\h\in\O$ we have
$$
\bbP(\s_{t,x}^\h\neq \s_{t,x}^+) = \bbP(\s_{t,x}^+=+1) -
\bbP(\s_{t,x}^\h=+1) = \frac12\, \big(\r_{t,x}(+)-\r_{t,x}(\h)\big)\,.
$$
Let $f$ be a function on $\O$ depending only on the spins in a finite
set $A\sset \Tree^b$ and let $\h\in \O_\a$. Then, using the invariance
of $\mu^+$, \ie $\mu^+P_t(f)=\int d\mu^+(\xi)(P_tf)(\xi)=\mu^+(f)$, 
for all $t$ large
enough depending on $A$, we have
\begin{gather*}
 |\bbE(f(\s^\h_t)) - \mu^+(f)|\leq | \bbE(f(\s^\h_t) - f(\s^+_t))|
 + |\int d\mu^+(\xi) \bbE(f(\s^+_t) - f(\s^\xi_t))| \cr 
  \leq 2\,\ninf{f}\sum_{x\in A}\Bigl[\ \bbP(\s_t^\h(x)\neq \s_t^+(x)) +
\int d\mu^+(\xi) \bbP(\s_t^\xi(x)\neq \s_t^+(x)) \ \Bigr] \cr
= \ninf{f}\sum_{x\in A}\Bigl[\ 2\,\r_{t,x}(+)-\rho_{t,x}(\h) -
 \mu^+(\s_x)\ \Bigr] %\cr
\leq \ninf{f} |A|\bigl[\  2 \,\nep{-\d t} + \nep{-t^\a}\ \bigr]\,. 
\end{gather*}
Therefore $\bbE(f(\s^\h_t)) \to \mu^+(f)$ for every bounded local
function and the weak convergence $\d_\h P_t\to\mu^+$ follows.
\end{proof}

Finally, the following generalization of Lemma  \ref{pro+} will
also be useful. Let us define the set $\O_{1,\d}$, for $\d>0$, as the set of $\h\in\O$ such that
(\ref{1}) above holds with the stretched exponential $\exp(-t^\a)$ replaced by the
true exponential $\exp(-\d t )$. Lemma \ref{pro+} then says that 
$+\in\O_{1,\d}$ for some $\d>0$.
\begin{lemma}
\la{mu+le}
For every $b,\b,h$, there exists $\d>0$ such that $\nu(\O_{1,\d})=1$ for any $\nu\geq\mu^+$.
\end{lemma}
\proof
Since $\nu\geq \mu^+$ we have $\nu(\r_{t,x})\geq \mu^+(\si_x)$ for all
$t\geq 0$. From Lemma \ref{pro+} we then infer
%Using monotonicity and with the constant $\d>0$ from Lemma \ref{pro+}
%we obtain %appearing in (\ref{2}), 
\begin{gather*}
\nu(|\rho_{t,x} - \mu^+(\s_x)|\geq \nep{-\d t/4})\leq   
\nep{\d t/2}\,\nu\bigl(|\rho_{t,x} - \mu^+(\s_x)|^2\bigr) \\
\leq \nep{\d t/2}\,\nu\bigl(\rho_{t,x}^2 - \mu^+(\s_x)^2\bigr)\,
\leq\, 2\,\nep{\d t/2}\,\bigl(\rho_{t,x}(+) -
  \mu^+(\s_x)\bigr) \leq 2\,\nep{-\d t/2}\,.
\end{gather*}
Therefore, the Borel-Cantelli lemma implies that
there exists $\e=\e(\d,b)>0$ and a subset $\O_0\sset \O$ of $\nu$-full
measure such that for all $\h\in \O_0$, all integers $j$ large enough
and all $x\in T_{\e j}$ (the tree of depth $\ell=\e j$), 
\begin{equation}
  |\rho_{j,x}(\h)- \mu^+(\s_x)| \leq \nep{-\d j/4}
\label{auxiliary}
\end{equation}
To prove the lemma we will 
establish a bound of the type (\ref{auxiliary}) on 
$|\rho_{t,r}(\h)-\mu^+(\s_r)|$, i.e.\ at the root
$x=r$, but for \emph{all} times $t$ large enough and not just integer
ones.
The case of general $x$ is obtained by straightforward modifications. 
%The result for arbitrary $x$ is obtained in the same way.
%a In order to get a similar bound on $|\rho_t^\h-\mu^+(\s_r)|$ for
%\emph{all} times $t$ large enough and not just integer ones, w
%
We simply write $\r_t$ for $\r_{t,r}$. Then, if $\inte{t}$ is the
integer part of $t$:
\begin{gather}
\rho_t(\h)-\mu^+(\s_r)= \rho_{\inte{t}}(\h)-\mu^+(\s_r) +
\int_{\inte{t}}^t
ds \, P_s g(\h)\,,\qquad g:=\cL\s_r %\\
%=\rho_{\inte{t}}(\h)-\mu^+(\s_r) + \int_{\inte{t}}^t
%ds \,\bbE(g(\s_s^\h)),\qquad g:=\cL\s_r
\label{probes}
\end{gather}
%If we let $G_s(\xi):= \bbE(g(\s_s^\xi))$, t
For $s\geq \inte{t}$ the 
Markov property yields $P_s g(\h) = P_{\inte{t}}P_{s-\inte{t}}g(\h)
= \bbE ([P_{s-\inte{t}} g] (\si_{\inte{t}}^\h))$.
%$\bbE(g(\s_s^\h))= \bbE\bigl(G_{s-\inte{t}}(\s_{\inte{t}}^\h)\bigr)$
On the other hand 
standard arguments (the so--called ``finite speed of propagation''
estimate) based on tail 
estimates for the mean one Poisson process (see e.g.  \cite{Mar}) show that
\begin{equation}
  \label{eq:G}
\sup_{0\leq u\leq 1} |P_ug(\xi)-P_ug(\xi')|\leq C_1\sum_x |\xi_x-\xi'_x|\,\nep{-C_2 \,d(x)}
\end{equation}
for constants $C_1,C_2$ with the property that
we can take $C_2$ as large as we wish provided $C_1$ is
large accordingly. Therefore
by the global coupling and the invariance of $\mu^+$, 
\begin{align*}
| P_s g(\eta) - \mu^+(g) |
& = | \int d\mu^+(\xi)
\bbE([P_{s-\inte{t}} g]
(\si_{\inte{t}}^\h) - [P_{s-\inte{t}}g](\si_{\inte{t}}^\xi) )|\\
& \leq C_1 \sum_x 
\,\nep{-C_2 \,d(x)}
\int d\mu^+(\xi)
\bbE(|\si_{\inte{t},x}^\h- \si_{\inte{t},x}^\xi|)\,.
\end{align*}
To handle the last term we add and subtract $\si_{\inte{t},x}^+$
so that by monotonicity
\begin{align*}
\int d\mu^+(\xi)
\bbE(&|\si_{\inte{t},x}^\h- \si_{\inte{t},x}^\xi|)
\leq 2\bbP(\si_{\inte{t},x}^+=1)-\bbP(\si_{\inte{t},x}^\h=1)
- \int d\mu^+(\xi)\bbP(\si_{\inte{t},x}^\xi=1)\\
& \leq  \frac12 \,|\r_{\inte{t},x}(\h) - \mu^+(\si_x))| 
+ \,(\r_{\inte{t},x}(+) - \mu^+(\si_x))\,.
\end{align*} 
Fix $j=\inte{t}$. When $x\in T_{\e j}$ we use (\ref{auxiliary})
for the first term above. If $\e$ is sufficiently small the argument
of Lemma \ref{pro+} also yields 
$$
 \r_{\inte{t},x}(+) - \mu^+(\si_x) \leq \nep{-\d_1\,t}\,,
$$
for some $\d_1>0$, uniformly in $x\in T_{\e \inte{t}}$. In conclusion,
for a suitable $\d_2>0$ we have 
\be
| P_s g(\eta) - \mu^+(g) |\leq \nep{-\d_2\,t} 
\la{probe}
\end{equation}
for all sufficiently large $t$.
The desired estimate now follows from (\ref{probe}) and (\ref{probes}), 
since $\mu^+(g)=0$ by invariance
of $\mu^+$. \qed

\section{Proof of Theorem \ref{theotree}}
\label{main_proof}
%Back to the proof of the theorem w
We will provide 
a unified proof of the three statements in Theorem \ref{theotree}.
In order to be able to do so we need some preliminary observations.
The first is that by the monotonicity of the events $\O_\a$ 
(Corollary \ref{incr}), statement
b) in 
Theorem \ref{theotree} is equivalent to 
\smallno 
{\em
\begin{enumerate}[b*)]
\item{For every $p>\frac12$ there exist
$b_0$ and $\b_0$ such that 
for $b\geq b_0$, $\b\geq \b_0$ and $h=0$, we have
$\bbP_p(\O_\a)=1$ for some $\a=\a(\b,b)>0$.}
\end{enumerate}
}
\medno 
Similarly, performing a global spin--flip, 
statement c) in Theorem \ref{theotree} can be rephrased as
\medno 
{\em
\begin{enumerate}[c*)]
\item For every $p>0$ there exist
$b_0$ and $\b_0$ such that 
for $b\geq b_0$, $\b\geq \b_0$ and $h=+h_c(\b,b)$, we have
$\bbP_p(\O_\a)=1$ for some $\a=\a(\b,b)>0$.
\end{enumerate}
}
\medno 
The last observation is that we may replace
statement a) in Theorem \ref{theotree} with \medno {\em
\begin{enumerate}[a*)]
\item For every $a>0$, $b\geq 2$,
there exist $p<1$ and $\b_0>0$ such that
for all $\b\geq \b_0$ and $h\geq -h_c(\b,b) + a$ 
we have $\bbP_p(\O_\a)=1$, for some $\a=\a(\b,h,b)>0$.

\end{enumerate}
}\medno 
In other words, we are taking $\b$ large enough.
To see why this is not restrictive recall that by an obvious 
domination argument one has $\bbP_p\geq \mu^+$ if
$$
p\geq p_{\b,h}:=\frac{\nep{(b+1+h)\b}}{\nep{(b+1+h)\b}+\nep{-(b+1+h)\b}},
\quad \text{\ie}\quad \b\leq \frac{1}{2(b+1+h)}\log(\frac{p}{1-p})\,.
$$
Lemma \ref{mu+le} therefore implies that $\bbP_p(\O_\a)=1$ for all 
$\a<1$ if $p\geq p_{\b,h}$. 
We then achieve the result of Theorem \ref{theotree} a) 
from a*) above by 
a suitable tuning of the parameter $p$.

%\medno  
%The above little discussion shows that, by possibly taking $p$ very
%close to one, it is sufficient to consider only
%the case of $h=0$ and $\b$ large enough.

\subsection{Main argument}
\la{main}
As the convergence result of Lemma \ref{pro+} makes clear, 
to prove Theorem \ref{theotree} we 
need to lower bound the quantity $\r_{t,x}(\h) -\mu^+(\s_x)$
in the three statements 
a*),b*) and c*) emphasized above. We shall focus only on the case $x=r$, 
since the case of arbitrary $x$ 
is obtained with essentially no modification. Setting
$\r_t(\eta):=\r_{t,r}(\eta)$ what we want is a bound of the form
\be
\r_t (\h) - \mu^+(\si_r) \geq -\,\nep{-t^\a}\,,
\la{toprobe}
\end{equation}
for all $t\geq t_0(\h)$, $\bbP_p$--almost all $\h$. 
As far as this section
goes, we shall not
distinguish the specific setting (a*,b* or c*), since all we do here
works for the 
three cases without any difference. What does depend on the setting 
are some key estimates that will be proved in the next two sections. 
The latter have been
emphasized as separate claims in the text (see Claims 1 to 4 below).

In order to describe the main idea behind 
the lower bound (\ref{toprobe}), we must first
introduce the notion of the Ising model and the associated Glauber
dynamics in a \emph{random environment of obstacles}.  
Realizations of the environment are described by elements $\o$ of
$\O$. We say that a vertex $x\in\bbT^b$ is an {\em
  obstacle} if $\o_x=-1$, and that $x$ is {\em free} if $\o_x=+1$.  We
call $T(\o)$ the largest connected component of the set of free vertices
containing the root. Note that $T(\o)=\emptyset$ if the root $r$ is
itself an obstacle. By construction, all vertices in $\partial T(\o)$
are obstacles. We will be mostly concerned with the case where $\o$ is
picked according to the product Bernoulli measure $\bbP_p$, i.e.\ when
each vertex is free with probability $p$, independently of all others.
In this case, $\bbP_p(T(\o) {\rm \;\,is\;\, infinite})$ is positive as
soon as $p>1/b$ and tends to $1$ as $p\nearrow 1$ for fixed $b$, 
or as $b\nearrow \infty$ for fixed $p$, see e.g.\ 
\cite{Peres}.

Given a realization of obstacles $\o$, the Ising model among obstacles
is defined as before by replacing the tree $\Tree^b$ with the random
tree $T(\o)$ and the configuration space $\O$ with the space
$$
\cB_\o:=\{\tau\in \O:\,
\tau_x=-1\,, \;\,\,\forall \,x\notin T(\o)\}\,.
$$
Given a finite subset $A\subset \Tree^b$ and $\tau\in\cB_\o$ we
denote by $\mu_{A,\o}^\tau$ the Gibbs measure $\mu_{A\cap T(\o)}^\t$.
We also write $\mu_{\ell,\o}^\tau$ for the Gibbs measure 
$\mu_{T_\ell(\o)}^\t$, where we use the notation $T_\ell(\o):= T_\ell \cap T(\o)$.
From this definition we see that obstacles act as a ``minus'' boundary
condition. The maximal allowed configuration $\tau\in
\cB_\o$ is such that $\tau=+1$ in $T(\o)$ and, when no confusion arises,
it will be always denoted by ``+''.  %If $T(\o)$ is infinite w
We will write
$\mu_\o^+$ for the Gibbs measure obtained as weak limit of
$\mu_{\ell,\o}^+$ as $\ell\to\infty$ (this is a finite volume Gibbs
measure with $-$ b.c.\ 
if $T(\o)$ is finite).  
%The above construction clearly
%makes sense also when $T(\o)$ is finite (but in this case $\mu^+_\o$ is
%nothing but $\mu_{T(\o)}^-$).

 \begin{figure}
  \centering
  \includegraphics[width=3in]{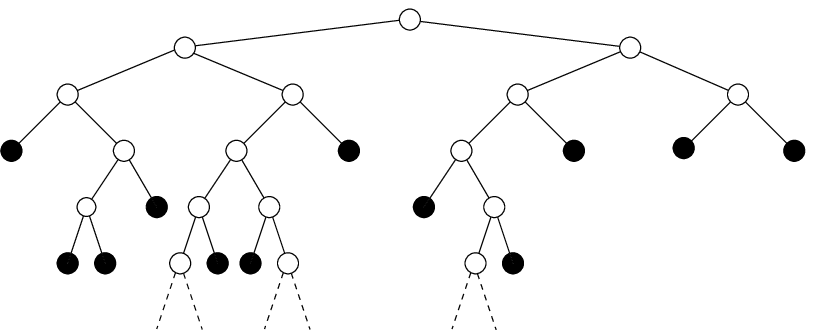}
  \caption{Free vertices ($\circ$) and obstacles ($\bullet$)
 in a given realization $\omega$ on the binary tree.} 
 \label{fig2}
  \end{figure}

%% %%%%%%%%%%%    PICTURE !@!@!!   %%%%%%%
%% \begin{figure}[h]%t= top, b = bottom, h = here
%% \centerline{
%% \psfig{file=albero.eps,height=1in}
%% }
%% \caption{Free vertices ($\circ$) and obstacles ($\bullet$)
%% in a given realization $\omega$ on the binary tree.} 
%% \label{fig2}
%% \end{figure}                    
%% %%%%%%%%%%    PICTURE !@!@!!   %%%%%%%%

Similar notations apply to the Glauber dynamics. 
Given a realization of
obstacles $\o$, $A\sset \Tree^b$ and $\t\in\cB_\o$, we will write
$\s_{t,\o}^{\xi,A,\t}$ for the Glauber dynamics in $A\cap T(\o)$ with
boundary condition $\t$ started from the restriction to $T(\o)$ of the
configuration $\xi\in\O$. If $A=T_\ell$ we simply write
$\s_{t,\o}^{\xi,\ell,\t}$. 
When $A=T(\o)$ the boundary condition is
necessarily $``-''$ and we will only write
$\s_{t,\o}^\xi$. In the same way we will use $\r_{t,\o}(\xi)$ for 
the expected value of $\s_{t,\o}^\xi$ at the root and
$\r_{t,\o}^{\ell,\t}(\xi)$ for the expected value of
$\s_{t,\o}^{\xi,\ell,\t}$
at the root.
Monotonicity implies that, for any
$\xi\leq \xi'$, $\s_{t,\o}^\xi \le
\s_{t}^{\xi'}$. In particular, 
\begin{equation}
  \label{eq:mono4}
\r_{t,\o}(+)\leq \r_{t}(\o)\,, \qquad t\geq 0\,.
\end{equation}

\medno 
We now turn to our main argument. 
Fix a length scale $\ell$, to be related later on to the time $t$,
a configuration $\eta\in \O$ and define the associated realization
$\o=\o(\h,\ell)$ of obstacles by the rule:
\begin{equation}
  \o_x = 
  \begin{cases}
  +1 & \text{ if $d(x)\leq \ell$} \\
  \h_x & \text{ otherwise }  
  \end{cases}
\la{rule}
\end{equation}
Clearly $\s_{t,\o}^\h \leq \s^\h_t$ so that
\begin{equation}
\r_t(\h) -\mu^+(\s_r)\geq [\r_{t,\o}(\h) -\mu^+_\o(\s_r)]
-[\mu^+(\s_r)-\mu^+_\o(\s_r)]\,.  
\label{main.01}
\end{equation}
If now $L=\ell^\g$, $\g>1$, is another length scale, monotonicity shows that
if we impose $+$ b.c.\ on the leaves of $T_L(\o)$ we may estimate
\begin{gather}
[\r_{t,\o}(\h) -\mu^+_\o(\s_r)]
-[\mu^+(\s_r)-\mu^+_\o(\s_r)] \nonumber \\ \ge
[\r_{t,\o}^{L,+}(\h) -\mu^+_{L,\o}(\s_r)] -
[\r_{t,\o}^{L,+}(\h)-\r_{t,\o}(\h)] 
-[\mu^+(\s_r)-\mu^+_\o(\s_r)]\,.
  \label{main.1}
\end{gather}
Notice that in the
above formula the role of the Bernoulli configuration $\h$
is twofold: it enters as the starting configuration in the first two
terms but it also defines the random realization of obstacles $\o$.

%                             %%%%%%%%%%    PICTURE !@!@!!   %%%%%%%
\begin{figure}[h]%t= top, b = bottom, h = here
\centerline{
\psfrag{ell}{$\ell$}
\psfrag{ELL}{$L$}
\psfrag{o}{$+\;+\;+\;+\;+\;+\;+$}
\psfig{file=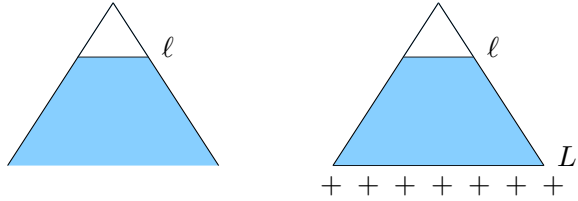,height=1in,width=3in}
}
\caption{Random obstacles below level $\ell$: 
Infinite tree (left) and
finite tree with $+$ boundary condition below
level $L=\ell^\g$, $\g>1$ (right).} 
\label{fig30}
\end{figure}                    %%%%%%%%%%    PICTURE !@!@!!   %%%%%%%
Most of the statements that will
be proved below on the r.h.s of (\ref{main.1}) concern properties which
hold almost surely with respect to the starting configuration $\h$ (and
therefore w.r.t. $\o$) picked according to the Bernoulli measure
$\bbP_p$. To simplify
the exposition, we shall adopt the following convention:
given some statements $\cE_\ell$, $\ell\in\bbN$, 
we say that $\cE_\ell$ holds $\bbP_p$--a.s.\ for $\ell$ sufficiently
large whenever $\eta\in\cE_\ell$ for all $\ell\geq\ell_0(\eta)$, for
$\bbP_p$--a.a. $\eta\in\O$, or in other words, $\bbP_p(\cE_\ell\; \text{eventually})=1$. 
We are now in a position to explain how we will bound the three
terms in the r.h.s of (\ref{main.1}).

%\medno
\smallskip

\subsubsection{{\bf Estimate on $[\mu^+(\s_r)-\mu^+_\o(\s_r)]$}}
Bounding the third term in (\ref{main.1})
is a purely static problem which on the tree can be
solved via a suitable recursion. In section \ref{claims} we prove 
\begin{claim}
  \begin{equation}
|\mu^+(\s_r)-\mu^+_\o(\s_r)|\leq \nep{-2\ell}\,,
\label{claimuno}  \end{equation}
\emph{$\bbP_p$--a.s.\ for $\ell$ 
sufficiently large.}
\end{claim}

\smallskip

\subsubsection{{\bf Estimate on $[\r_{t,\o}^{L,+}(\h)
    -\mu^+_{L,\o}(\s_r)]$}}
The first term in (\ref{main.1}) is related to the speed of relaxation to
  equilibrium in the finite tree $T_L(\o)$ with plus b.c.
Here we need the following bound
 on the logarithmic Sobolev constant
  $\csob(\mu_{L,\o}^+)$. 
{\em
  \begin{claim}
There exists $\zeta<\infty$ independent of $\ell$ such that 
\begin{equation*}
\csob(\mu_{L,\o}^+)\geq L^{-\zeta}
\end{equation*}
holds $\bbP_p$--a.s.\ for any $L\geq \ell$ sufficiently large.
\end{claim}
}
We can
repeat the computation (\ref{3}) with $\mu_\ell^+$ replaced by
$\mu_{L,\o}^+$, $\phi_t^{\ell,+}$ by $
\phi_{t,\o}^{L,+}(\eta):=\r_{t,\o}^{L,+}(\eta) -\mu^+_{L,\o}(\s_r)$ 
and $q_t := 1+\nep{2\,t\,\csob(\mu_{L,\o}^+)}$ to obtain
\be
|\phi_{t,\o}^{L,+}(\h)| %\leq (\mu_{\ell}^+(+))^{-\frac1{q_t}}
\leq \exp{(C_1b^L\nep{-\csob(\mu_{L,\o}^+) \,t})} \,
\|\phi_{t/2,\o}^{L,+} \|_{2,\mu_{L,\o}^+}\,.
\label{step1.1}
\end{equation}
for some constant $C_1<\infty$.
Assuming Claim 2 above, using $\cgap \geq 2\csob$ (which is always
true, see e.g.\ \cite{Saloff}) we estimate (\ref{step1.1}) with the
help (\ref{gapvar}) and obtain, for $L=\ell^\g$,
\be
|\r_{t,\o}^{L,+}(\h)
    -\mu^+_{L,\o}(\s_r)| = |\phi_{t,\o}^{L,+}(\h)|
\leq \exp{(C_1b^L\nep{-t/\ell^{\g \,\zeta}})}
\nep{-t/\ell^{\g \,\zeta}}
\la{eq:main_first}
\end{equation}
$\bbP_p$--a.s.\ for $\ell$ sufficiently large. 

\subsubsection{{\bf Estimate on $[\r_{t,\o}^{L,+}(\h)-\r_{t,\o}(\h)] $}}
The control of the second term in (\ref{main.1}) is a true
  dynamical question and it involves proving that the two processes
  $\s_{t,\o}^\h$ and $\s_{t,\o}^{\h,L,+}$ remain identical at the root
  up to time $t$ with large probability.  This is  achieved via a
  coupling argument together with some equilibrium estimates. The final
  bound will be of the form
  \begin{equation}
    \label{eq:main_second}
    \r_{t,\o}^{L,+}(\h)-\r_{t,\o}(\h) \leq t\,\nep{-2\ell}
  \end{equation}
 $\bbP_p$--a.s.\ for $\ell$ sufficiently large, 
provided that $\g>\zeta+1$, where $\zeta$
is the constant appearing in (\ref{eq:main_first}).
The argument goes as follows.

To keep the notation to a minimum, we will
abbreviate the two processes $\s_{t,\o}^\h$ and $\s_{t,\o}^{\h,L,+}$ with 
$\xi_{t}^1$ and
$\xi_t^2$ respectively. Set 
\begin{align*}
\L &:=\{x\in T(\o):\ d(x)=2\ell\}\,,\quad 
\bar\L :=\{x\in T(\o):\   \frac 32\ell \leq d(x)\leq \frac 52\ell\}\,.
\end{align*}
By the global coupling %clearly
\begin{align}
0\leq  & = \r_{t,\o}^{L,+}(\h)-\r_{t,\o}(\h) =  \bbE
\left[\xi_{t,r}^2-\xi_{t,r}^1\right]\nonumber\\
& \leq
2\,\bbP\Big[\exists s\leq t\,,\; \exists \, x\in\L\,:\;\,
\xi^{1}_{s,x}\neq \xi^{2}_{s,x}\, 
\Big]\,.
 \label{7}
\end{align}
Define $A_j$, $j=1,2,\dots,\inte{t}$, as the event
$$
A_j = \{\exists x\in\bar \L\,:\;\, \xi^{1}_{j,x} \neq \xi^{2}_{j,x}
\,\}\,.
$$
The r.h.s.\ in (\ref{7}) is then estimated from above by
\begin{equation}
\sum_{j=1}^{\inte{t}} \bbP(A_j)  \,+\,
\sum_{j=1}^{\inte{t}+1} \bbP\Big[A_{j-1}^c\cap B_j\Big]
\label{8}
\end{equation}
where 
$$
B_j:=\{\exists s\in[j-1,j]\,,\;\exists x\in\L\,: \;
\xi^{1}_{s,x}\neq \xi^{2}_{s,x}\,\}\,.
$$
The probability of the event $A_{j-1}^c\cap B_j$ 
is estimated by a standard argument: the event $A_{j-1}^c\cap B_j$ implies that a
discrepancy between $\xi_{j-1}^{1}$ and
$\xi_{j-1}^{2}$ located
outside $\bar \L$, reaches in a time smaller than $1$ 
a point $x\in\L$. Since there are at most $b^{\frac 52\ell}$ possible 
(self--avoiding) paths 
from $\L$ to $(\bar \L)^c$ and since the rates are bounded by one, 
a simple tail estimate for Poisson random variables implies
\begin{equation}
\bbP(A_{j-1}^c\cap B_j)\leq c\,b^{\frac 52 \ell}\,
\nep{-\frac 12 \ell \log(\ell/ c)}\leq \nep{-3\ell}
\label{9}
\end{equation}
for a suitable constant $c$ and all sufficiently large $\ell$.  
Similarly, if we look at the event $A_j$,
we are requiring that at least one of the  discrepancies at time $0$ in level $L$ travels 
up to level $\frac 52 \ell$ in a time less than $j$. Therefore
\begin{equation}
\bbP(A_{j}) \leq c \, b^{L}\,\nep{-(L-\frac 52 \ell)\log((L-\frac 52
  \ell)/cj)}\leq \nep{-L}\,,\qquad \forall j  < \e L\,.
\label{11}
\end{equation}
for some $c<\infty$ and for all $\e=\e(b,c)$ sufficiently small.
We are therefore left with the  estimate of $\bbP(A_{j})$ 
for $j \geq \e L$. The argument for this case goes as follows.

Fix a point
$x\in\bar\L$ and recall that $x$ is at some level between $\frac 32
\ell$ and $\frac 52 \ell$. Let $r_x=r(x,\ell)$ be the ancestor of $x$ at
level $\ell$ and let $T_{r_x}(\o)$ be the subtree of $T(\o)$
rooted at $r_x$ and containing all descendants of $r_x$. 
Let also $T_{r_x,h}(\o)$, $h\in\bbN$, be the finite subtree of
$T_{r_x}(\o)$ obtained by considering only the first $h$ levels of 
$T_{r_x}(\o)$
(so that $T_{r_x,h}(\o)$ consists of $r_x$ and all its descendants in $T(\o)$
lying between level $\ell$ and $\ell+h$). 
When $h=2\ell$, $T_{r_x,2\ell}(\o)$ is the tree between levels $\ell$
and $3\ell$, so that $\bar\L\subset T_{r_x,2\ell}(\o)$ and 
$d(\bar\L,(T_{r_x,2\ell}(\o))^c)\geq \frac\ell 2$. 
%                             %%%%%%%%%%    PICTURE !@!@!!   %%%%%%%
\begin{figure}[h]%t= top, b = bottom, h = here
\centerline{
\psfrag{z}{$r_x$}
\psfrag{ell}{$\ell$}
\psfrag{ell1}{$\frac32\ell$}
\psfrag{ell2}{$2\ell$}
\psfrag{ell3}{$\frac52\ell$}
\psfrag{ell4}{$3\ell$}
%\psfrag{la}{$\La\;\rightarrow$}
\psfrag{.}{{\Large $\cdot$}} \psfrag{x}{$x$}
\psfig{file=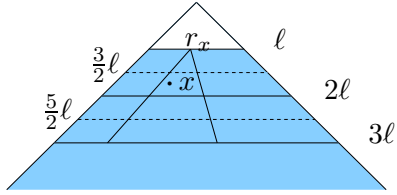,height=1in,width=2in} }
\caption{The vertex $x\in\bar \L$ and the associated tree $T_{r_x,2\ell}$.} 
\label{fig3}
\end{figure}                    %%%%%%%%%%    PICTURE !@!@!!   %%%%%%%
Call $\nu_{r_x,h}^{+,+}$ ($\nu_{r_x,h}^{-,+}$)
the Ising--Gibbs measure on $T_{r_x,h}(\o)$ with $+\,(-)$ b.c. 
above the root $r_x$ and $+$ b.c. below the leaves at level $\ell+h$ and
let $\nu_{r_x,\infty}^{-,+}=\lim_{h\to \infty} \nu_{r_x,h}^{-,+}$.

Finally, we denote by
$\xi^{3}_t$ the Glauber dynamics evolving in $T_{r_x,2\ell}(\o)$
with $+$ boundary conditions both above the root $r_x$ and below the leaves
of $T_{r_x,2\ell}(\o)$ and with initial configuration $\h$. Notice that in fact
$\xi^{3}_t$ starts from all pluses because, by construction,  $\h(y)=+1\ \forall y\in
T(\o), \; d(y)\geq \ell$.  

With the above notation and using monotonicity we can now write 
\begin{align}
\bbP\big[\xi^{1}_{j,x}
\neq\xi^{2}_{j,x} 
\big] & = \bbP\big[\xi^{2}_{j,x} = +1\big]
- \bbP\big[\xi^{1}_{j,x} = +1\big]
\nonumber\\
&\leq \bbP\big[\xi^{3}_{j,x} = +1\big]
- \nu_{r_x,\infty}^{-,+}(\si_x = +1)\,.
\label{13}
\end{align}
The r.h.s.\ in (\ref{13}) is then decomposed into the sum of two terms:
\begin{equation}
\bbP\big[\xi^{3}_{j,x} = +1\big]
- \nu_{r_x,2\ell}^{+,+}(\si_x = +1)\,,
\label{13_1}
\end{equation}
and
\begin{equation}
\nu_{r_x,2\ell}^{+,+}(\si_x = +1) - \nu_{r_x,\infty}^{-,+}(\si_x = +1)\,.
\label{13_3}
\end{equation}
In order to bound the term in (\ref{13_1}) we would like to argue as 
in (\ref{step1.1}) and therefore we need the following:
{\em
  \begin{claim}
There exists $\zeta<\infty$ 
such that 
\begin{equation*}
\min_{x\in \bar \L}\,\csob(\nu_{r_x,2\ell}^{+,+})\geq \ell^{-\zeta}
%\eqno({\bf C}\,3)
%\label{sobbound2}
\end{equation*}
holds $\bbP_p$--a.s.\ for $\ell$ sufficiently large.
  \end{claim}
}
The argument used in (\ref{step1.1}) now yields that the first term
(\ref{13_1}) satisfies
\begin{equation}
\bbP[\xi^{3}_{j,x} = +1]
- \nu_{r_x,2\ell}^{+,+}(\si_x = +1)
\leq  \exp{\big(c\, b^{3\ell}\nep{-j /\ell^{\zeta}}\big)}
\, \exp{\big(-j /\ell^{\zeta}\big)}\,.
\label{14}
\end{equation}
Therefore, if $\g > \zeta+1$, using
$j \geq \e \ell^\g $ we may write
\begin{equation}
\bbP(\xi^{3}_{j,x} = +1)
- \nu_{r_x,2\ell}^{+,+}(\si_x = +1) \leq b^{-3\ell} \nep{-3\ell}
\label{15}
\end{equation}
for $\ell$ large enough (independent of $x$).
In conclusion, for any
$\g>\zeta+1$, 
\begin{equation}
  \sum_{x\in \bar \Lambda}\left\{ \bbP[\xi^{3}_{j,x} = +1]
- \nu_{r_x,2\ell}^{+,+}(\si_x = +1) \right\}
\leq \nep{-3\ell}
\la{claim9}
\end{equation}
$\bbP_p$--a.s.\ for $\ell$ sufficiently large.

\smallskip

As far as the term (\ref{13_3}) is concerned we will establish:
{\em
  \begin{claim}
\begin{equation*}
\sum_{x\in \bar \Lambda}\left\{ \nu_{r_x,2\ell}^{+,+}(\si_x = +1) -
\nu_{r_x,\infty}^{-,+}(\si_x = +1)\right\} \leq \nep{-3\ell}\,,  
\end{equation*}
$\bbP_p$--a.s.\ for $\ell$ sufficiently large.
\end{claim}
}
Collecting (\ref{claim9}) and Claim 4, we have thus shown that 
$\bbP(A_j) \leq 2\nep{-3\ell}$, $j\geq \e L=\e \ell^\g$.
%
%$p_0\in(0,1)$, $\b_0\in(0,\infty)$, such that 
%for all $p\geq p_0$, $\b\geq \b_0$, if $\g>a(\b)+1$  
%\begin{equation*}
%  \bbP(A_j) \leq \nep{-\ell} \quad \forall j\geq \e \ell^\g
%\end{equation*}
%$\bbP_p$--a.s.\ for $\ell$ sufficiently large. 
Together 
with (\ref{9}) and (\ref{11}) this  completes the proof of
(\ref{eq:main_second}).

\subsubsection{{\bf Conclusion}}
In conclusion, from (\ref{main.1}), using 
(\ref{claimuno}), (\ref{eq:main_first}) and
(\ref{eq:main_second}) we have the bound
$$
\r_t (\h) - \mu^+(\si_r) \geq - (\nep{-2\ell} + 
\exp{(C_1b^L\nep{-t/\ell^{\g \,\zeta}})}
\nep{-t/\ell^{\g \,\zeta}}
+ t\,\nep{-2\ell})\,.
$$
If e.g.\ $\ell =t^{\a}$ with $\a>0$ such that 
$\a\g(1+ \zeta)<1$ and 
$\g>\zeta+1$, then 
\begin{equation}
  \label{eq:main_final}
  \r_t (\h) - \mu^+(\si_r) \geq -\nep{-t^\a } \,,
\end{equation}
$\bbP_p$--a.s.\ for $t$ sufficiently large. 
Therefore the three steps above are
sufficient to end the proof of Theorem \ref{theotree}. 
Note that the coefficient 
$\a$ depends on the various parameters ($b,\b,h$ etc.) only via the 
constant $\zeta$ coming from the logarithmic Sobolev inequality
in Claim 2.

%%%%%%%%%%%%%%%%%%%%%%%%%%%%%%%%%%%%%%%%%%%
%%%%%%%%%%%%%%%    OBSTACLE.TEX
%%%%%%%%%%%%%%%%%%%%%%%%%%%%%%%%%%%%%%%%%%%

\section{Recursive analysis among obstacles}
\label{recursive}
In this section we establish a number of key estimates for the
Ising Gibbs measure among obstacles. Once these results are
established it will be rather easy to prove Claims
1 to 4 (see next section).

As in the previous section $\o\in\O$ will denote a random realization of the 
obstacle--environment and $\mu_\o^+$ the associated Ising 
plus phase. We emphasize however that here, contrary
to (\ref{rule}), $\o$ is picked 
according to the product Bernoulli measure 
$\bbP_p$ on the \emph{whole tree $\Tree^b$}, \ie \emph{each} vertex $x\in
\Tree^b$ 
is free with probability $p$ independently of all others. 
% and occupied (\ie carries a minus boundary condition),
%
%
%Essentially two cases will be considered:
%either $b$ is given and $p$ is large, or $p$ is given and $b$ is large.
%Clearly, in both cases with positive $\bbP_p$--probability a large
%number of rays of the tree $T(\o)$ will reach infinity. Quantitative versions
%of such a statement will be used to derive the desired 
%properties for the equilibrium measure $\mu_\o^+$. 
%We shall not attempt to obtain sharp
%relations between the various parameters $p,b,\b,h$ and our arguments
%will be essentially perturbative.
 
\subsection{Coupling coefficients and path weights}
As in the homogeneous case treated in \cite{MaSiWe}, 
the analysis of equilibrium properties 
is reduced to the study of certain
coupling coefficients. \medno For a given $\o$ we define the ratio 
\begin{equation}
R(\o)=\frac{\mu_\o^+(\si_r=-1)}{\mu_\o^+(\si_r=+1)} \,.
\la{r1}
\end{equation}
We agree that $R(\o)=\infty$ if $\o_r=-1$. 
For every $z\in\bbT^b$ we set $R_z(\o):=R(\theta_z\o)$, where $\theta_z$
denotes the shift induced by the natural group action on the tree:
$(\theta_z\o)_x=\o_{z+x}$. If $\o_z=(\theta_z\o)_r = +1$ and 
$z_1,\dots,z_b$ denote the children of $z\in T(\o)$, 
one has the following easily checked recursive 
relation (see e.g.\ \cite{MaSiWe}):
\begin{equation}
R_z(\o) = %\nep{-2\b h}
\e^h
\prod_{k=1}^b F_\b(R_{z_k}(\o))\,,\quad 
\la{w3}
\end{equation}
where, from now on, we use the following notation
\be
F_\b(a)= \frac{\e+a}{1+\e a}\,,\quad\;\e:=\nep{-2\b } 
\la{defFe}
\end{equation}
%The above notation $\e:=\nep{-2\b }$ will be used 
%throughout the rest of the paper. 
To illustrate the use of the variable $R$ defined in (\ref{r1}), consider
a vertex $z\in T(\o)$ together with one of its ancestors $y$ and
denote by $\mu^{y,+}_{\o}$ (resp.\ $\mu^{y,-}_{\o}$) the %marginal
%at $z$ of 
measure $\mu_\o^+$ conditioned to have $\si_y=+1$ (resp.\ $\si_y=-1$).
Suppose we want to compute 
the total variation distance between the marginals at the vertex $z$, 
which we denote by $\|\mu^{y,+}_{\o}- \mu^{y,-}_{\o}\|_z$.
Since the spin at $z$ can take only two values, the latter equals 
$\mu^{y,+}_{\o}(\si_z=+1)-\mu^{y,-}_{\o}(\si_z=+1)$. 
If $y$ is the parent of $z$, using 
$\mu_{\o}^{y,\pm}(\si_z=+1) = (\e^{\pm} R_z(\o) + 1)^{-1}$
we see that 
\begin{equation}
%\mu_\o(\si_z=+1\tc\eta_y=+1)-\mu_\o(\si_z=+1\tc\eta_y=-1)
\|\mu^{y,+}_{\o}- \mu^{y,-}_{\o}\|_z = K_\b(R_z(\o))\,,\quad\; 
\la{totvar1}
\end{equation}
where the function $K_\b:[0,\infty)\to [0,1]$ is defined by 
\begin{equation}
K_\b(a)=\frac{1}{\e \,a + 1}\,-\,\frac{1}{\e^{-1}a + 1}\,.
\la{k1}
\end{equation}
For every $\ell\in\bbN$ we define 
the set of descendants of $y$ at depth $\ell$:
\begin{equation}
D_{y,\ell}(\o)= \{x\in T(\o)\;{\rm descendant\;of\;} y\,:\;d(y,x)=\ell\}
\la{dyl}
\end{equation}
To compute the total variation distance $\|\mu^{y,+}_{\o}-
\mu^{y,-}_{\o}\|_x$ for some $x\in D_{y,\ell}(\o)$ we may proceed as
follows.  Let $z_1,\dots,z_\ell=x$ be the vertices along the path from
$y$ to $x$. We couple the measures $\mu^{y,+}_{\o}, \mu^{y,-}_{\o}$
recursively in such a way that, for every $i<\ell$, given that the
corresponding configurations coincide at $z_i$ then they coincide at
$z_{i+1}$ with probability $1$, while given that there is disagreement
at $z_i$ then disagreement persists at $z_{i+1}$ with probability
$\|\mu^{z_i,+}_{\o}- \mu^{z_{i},-}_{\o}\|_{z_{i+1}}=
K_\b(R_{z_{i+1}}(\o))$. In this way the probability of a disagreement
percolating {\em down} the tree from $y$ to $x$ equals 
\begin{equation}
\|\mu^{y,+}_{\o}- \mu^{y,-}_{\o}\|_x = \prod_{i=1}^\ell K_\b(R_{z_i}(\o))
\la{totvar2}
\end{equation}
Moreover, if $|\si-\si'|_{y,\ell}$ denotes the Hamming distance 
(counting the number of disagreements) between
$\si$ and $\si'$ restricted to the set $D_{y,\ell}(\o)$,
%of descendants $x\in T(\o)$ of $y$ 
%such that $d(x,y)=\ell$, 
the above argument implies that we can find 
a coupling $\nu_\o$ of $\mu^{y,+}_{\o},\mu^{y,-}_{\o}$
such that the expected value of $|\si-\si'|_{y,\ell}$ satisfies
\begin{equation}
\nu_\o(|\si-\si'|_{y,\ell})\leq \sum_{x\in D_{y,\ell}} W_\o(\G_{y,x})\,,
\la{ham}
\end{equation}
where we introduced the path $\G_{y,x}$ between $y$ and $x$,
consisting of the sites $z_1,\dots,z_\ell=x$ as above, and the associated
weight
\begin{equation}
W(\G_{y,x},\o)= \prod_{i=1}^\ell K_\b(R_{z_i}(\o))\,.
\la{w1}
\end{equation}
The rest of this section is concerned with estimates showing that, 
in a suitable sense, $R(\o)$ 
and $W_\o$ are small with large probability.\subsection{Estimates on $R$}
%For every $z\in \bbT^b$ we call $\bbT_z$ the subtree 
%rooted at $z$, i.e.\ the subtree consisting of $z$ and all its descendants. 
%Note that, as random variables on the probability space 
%$(\O,\bbP_p)$, $R_z$ and $R_y$ are independent whenever 
%$\bbT_z\cap \bbT_y=\emptyset$. 
We write
$\wt\bbP_p$ for 
the probability $\bbP_p$ conditioned to have $\o_r=+1$.
We want an estimate of the type 
\begin{equation}
\wt\bbP_p \left(R \geq \e \right)\leq \d\,,
\la{R1}
\end{equation}
where $\e=\nep{-2\b}$ and $\d$ is a small parameter. We start with the
setting of statement \emph{a*} in the proof of Theorem 
\ref{theotree}.
\begin{lemma}
\la{R_a}
For any $\d>0$, $a>0$, $b\geq 2$, 
there exist $p_0<1$ and $\b_0<\infty$ such that
(\ref{R1}) holds for all $p\geq p_0$, $\b\geq \b_0$ and $h\geq -h_c(\b)+a$.
\end{lemma}
\proof
For any integer $\ell$ we define
\begin{equation}
R^{\ell}(\o) = \frac{\mu_{\ell,\o}^+(\si_r=-1)}
{\mu_{\ell,\o}^+(\si_r=+1)} \,.
\la{r2}
\end{equation}
Since $\mu_{\ell,\o}^+\to\mu_\o^+$, 
we have $R^{\ell}\to R$, $\ell\to\infty$, $\bbP_p$--a.s.
Moreover, monotonicity implies 
$R^\ell(\o)\leq R^{\ell+1}(\o)$, so that the convergence is monotone. 
Then it is sufficient to establish (\ref{R1}) for 
$R^{\ell}$ in place
of $R$, uniformly in $\ell$. We will give the proof 
only in the case $b=2$, since all the estimates below are easily adapted 
to the case of larger values of $b$. Recall that in general 
(see e.g.\ \cite{Georgii}) one has $h_c(\b)=(b-1)+O(\b^{-1})$,
so that, replacing $a$ with $2a$ and taking $\b$ sufficiently large,
%for the rest of this proof 
we can assume $h\geq -1 + a$ without loss of generality.  
Let us define the probabilities
\begin{equation}
q^{(k)}_\ell=\wt\bbP_p\left(R^{\ell}>2^{-2k}\e^{1-ka}\right)\,,
\quad k=0,1,2,\dots
\la{qk}
\end{equation}
Let now $z_1$, $z_2$ denote the two children of the root
and observe that the corresponding ratios 
$R_i$, $i=1,2$ are i.i.d.\ random variables with the same distribution
as $R^{\ell-1}$. On the event $\{\o_r=+1\}$ the basic relation (\ref{w3})
applies and we have
\begin{equation}
R^\ell \leq \e^{a-1}\,F_\b(R_1)\,F_\b(R_2)\,.
\la{rel}
\end{equation}
Using the uniform bound $F_\b\leq \e^{-1}$, 
we see that, in particular, $R^\ell\leq \e^{a-1}\e^{-2}$. Therefore 
$q^{(k)}_\ell = 0$ for every $\ell\geq 1$
as soon as $k>k_0:=\inte{\frac4a -1}$ and $\b$ is large enough. 
Suppose now $R_1\leq \e$. Then by (\ref{rel}),
using $F_\b(R_1)\leq 2\e$  we have 
$$
R^\ell\leq 2\,\e^{a}F_\b(R_2)\,.
$$
Since $F_\b(t)\leq \e + t$, the event $R^\ell\geq \e$ forces 
$R_2\geq \frac{\e^{1-a}}4$ for large enough $\b$. 
Considering also 
the event $\{\o_{z_1}=-1\}\cup
\{\o_{z_2}=-1\}$ and the event 
$\{\o_{z_1}=+1,R_1> \e\}\cap \{\o_{z_2}=+1,R_2>\e\}$ we may then 
use (\ref{qk}) to write
\begin{equation*}
q^{(0)}_\ell\leq 2(1-p) + (q^{(0)}_{\ell-1})^2 + 2\,q^{(1)}_{\ell-1}\,.
\la{qq1}
\end{equation*}
The same reasoning as above actually shows that for any $k$ one has
\begin{equation}
q^{(k)}_\ell\leq 2(1-p) + (q^{(0)}_{\ell-1})^2 + 2\,q^{(k+1)}_{\ell-1}\,.
\la{qq2}
\end{equation}
%In particular, since $q^{(k_0+1)}_\ell=0$, we have 
%\be
%q^{(k_0)}_\ell\leq 2(1-p) + (q^{(0)}_{\ell-1})^2\,, 
%\la{qq3}
%\end{equation}
%for every $\ell$.
From the monotonicity in $\ell$ of $R^\ell$ we see that $q^{(k)}_{\ell-1}\leq
q^{(k)}_{\ell}$ for any $k$ and $\ell$. Therefore a simple iteration of 
(\ref{qq2}) gives that 
$$
q^{(0)}_\ell\leq \sum_{m=0}^{j-1}2^m
\left\{2(1-p) + (q^{(0)}_{\ell-1})^2\right\} + 2^j\,q^{(j)}_{\ell-1}\,,
$$
for any $j=1,2\dots$. When $j=k_0+1$,
$q^{(j)}_{\ell-1}=0$ and we have the recursive estimate
\begin{equation}q^{(0)}_\ell\leq 
2^{k_0+1}
\left\{2(1-p) + (q^{(0)}_{\ell-1})^2\right\}\,.
\la{qq4}
\end{equation}
This implies that for every $\d>0$ we can choose $p_0<1$ and $\b_0<\infty$
such that $q^{(0)}_\ell\leq \d$, for every $\ell\geq 1$, $p\geq p_0$
and $\b\geq \b_0$. To see this simply observe that when $\ell=1$
the ``$+$'' boundary condition imposes 
$R_{i}=0$ on every child $z_i$ such that
$\o_{z_i}=+1$ and therefore $q^{(0)}_{1}\leq 2(1-p)$, which can be 
made arbitrarily small. Thus, climbing up the tree with the relation 
(\ref{qq4}), we see that $\sup_{\ell\geq 1}q^{(0)}_\ell\leq \d$
as soon as e.g.\ $p\geq 1-1/(2^{5+2k_0})$. \qed

\bigskip 

We turn to the setting of statement
\emph{b*}
in the proof of Theorem \ref{theotree}.
\begin{lemma}
\la{R_b}
For any $\d>0$, $p>\frac12$, there exist $b_0\in\bbN$, $\b_0<\infty$ and $c>0$ such that
(\ref{R1}) holds for all $\b\geq \b_0$, $b\geq b_0$ and $h=0$, with $\d=\nep{-c b}$. 
\end{lemma}
\proof
Recall the definition (\ref{r2}) of $R^\ell(\o)$. Let $z_1,\dots,z_b$
denote the children of the root $r$ and let $m(\o)$ stand for the number 
of obstacles among them: $m(\o)=\sum_{i=1}^b\un_{\{\o_{z_i}=-1\}}$.
Since $p>\frac12$, from standard large deviation estimates for
the binomial distribution there exist positive numbers $a_1,a_2>0$
and $b_0\in\bbN$
such that 
such that 
\begin{equation}
\bbP_p\left( m \geq \left(\frac12-a_1\right)b \right)\leq \nep{-a_2 b}\la
{expb}\,,
\end{equation}
for all $b\geq b_0$.
Suppose now that $\o_r=+1$ and 
$m < (\frac12-a_1)b$, \ie the root has at least 
$\frac b2 + a_1b$ free children. Suppose only one of these
free children, say $z$, is such that the associated ratio $R_z$
satisfies $R_z\geq \e$. In this case 
(\ref{w3}) yields 
\begin{equation}
R\leq \e^{-m} (2\e)^{b-m-1} F(R_z)\leq 2^b \e^{2a_1b-1} F_\b(R_z)\,.
\la{qb1}
\end{equation}
Since $F_\b\leq\e^{-1}$ it is clear that we can take $b_0,\b_0$ so large
that in the above situation it is impossible to have $R\geq \e$ for
all $b\geq b_0$ and $\b\geq \b_0$. The above discussion says,
in particular, that if $m < (\frac12-a_1)b$ and $R\geq \e$, 
then there must be at least $2$ children of $r$ with ratio smaller than
$\e$. Thus, recalling the definition of the probabilities $q^{(0)}_\ell$
(\ref{qk}), we obtain
\begin{equation}
q^{(0)}_\ell
\leq \nep{-a_2 b} + \sum_{n=2}^b \binom{b}n \left(q^{(0)}_{\ell-1}\right)^n\,.
\la{qb}
\end{equation}
Because of the $+$ boundary condition at level $\ell$,
the argument of (\ref{qb1}) gives $q^{(0)}_{1}\leq \nep{-a_2 b}$.
The claim then follows by induction: Suppose 
$q^{(0)}_{\ell -1}\leq \d$ with $\d:=\nep{-a_2 b/2}$. 
Then (\ref{qb}) implies $q^{(0)}_\ell
\leq \nep{-a_2 b} + (1+\d)^b - (1+\d b) \leq \nep{-a_2 b} + \frac12 \,b^2\d^2$,
and therefore $q^{(0)}_\ell\leq \d$ for 
$b$ suitably large. 
\qed

\bigskip

Finally, for the statement
\emph{c*} in the proof of Theorem \ref{theotree} we need the following
\begin{lemma}
\la{R_c}
For any $\d>0$, $p> 0$, 
there exist $b_0\in\bbN$, $\b_0<\infty$ and $c>0$ such that
(\ref{R1}) holds for all $\b\geq \b_0$, $b\geq b_0$ and $h=+h_c(\b)$,
with $\d=\nep{-c b}$. 
\end{lemma}
\proof
As in the previous proof we denote by $m(\o)$ the number of
obstacles among the children of the root. Since $p>0$, 
%from standard large deviation estimates for
%the binomial distribution 
there exist positive numbers $a_1,a_2>0$ and $b_0\in\bbN$
such that 
\begin{equation}
\bbP\left( m \geq \left(1-a_1\right)b \right)\leq \nep{-a_2 b}
\la{expc}\,,
\end{equation}
for all $b\geq b_0$.
We recall that $h_c(\b)=(b-1)+O(\b^{-1})$. In particular, we may assume 
without loss of generality that 
the magnetic field satisfies $h\geq b-2$. If $z_1,\dots,z_b$ denote the 
children of the root, by (\ref{w3}) we then have
$$
R\leq \e^{b-2}F_\b(R_{z_1})\cdots F_\b(R_{z_b})\,.
$$
Reasoning as in (\ref{qb1}) we see that
if $\o_r=+1$ and $m(\o) < (1-a_1)b$, then we must have 
more than one free children with ratio smaller than $\e$ in order to produce
the event $R\geq \e$. It follows that we may estimate 
the probabilities $q^{(0)}_\ell$ exactly as in (\ref{qb}).
When $\ell=1$ the $+$ boundary condition implies 
$R\leq \e^{b-2}\e^{-m}\e^{b-m}$. Therefore on the event $m(\o) < (1-a_1)b$
it is impossible (for suitably large $b,\b$) to have $R\geq \e$.
This gives $q^{(0)}_{1}\leq \nep{-a_2 b}$. As in the proof of Lemma \ref{R_b}, 
the desired result now follows by induction.
\qed

%%%%%%%%%%%%%%%%%%%%%%%%%%%%%%%%%%%%%%%%%%%%%%%%%%%%%%
\subsection{Estimates on $W$}
We turn to an estimate on the weight $W$ introduced in (\ref{w1}). 
Recall the definition (\ref{dyl}) of the set $D_{y,\ell}(\o)$. 
Below we simply write $D_\ell = D_\ell(\o)$ when $y$ coincides with the root $r$.
We also write $W(x):=W(\G_{r,x},\o)$,
for any $x\in \bbT^b$. 
%If $x\notin T(\o)$
%then we set $W(x)=0$. 
We look for an estimate of the form: There exists $t_0>0$ such that 
\begin{equation}
\bbE_p \left[
\,\exp\,\left( t\sum\nolimits_{x\in D_\ell } W(x) \right)\right]
\leq 2\,,%\frac1{1-2t}
\la{expw}
\end{equation}
for every $t\leq t_0$, $\ell\geq 1$. The value of $t_0$ will depend on
the parameters $a,b,\b$ in case {\em a*}, on $b,\b$ and $p$ in cases
\emph{b*,c*}.

\smallskip

We start with the setting of statement \emph{a*}.
\begin{lemma} \la{weights_a} For any $a>0$, $b\geq 2$, there exist
  $\b_0<\infty$ and $p_0<1$, such that (\ref{expw}) holds for any
  $\b\geq \b_0$, $p\geq p_0$, $h\geq -h_c(\b)+a$ .
\end{lemma}
\proof 
The main difficulty in proving (\ref{expw}) is the non--independence of the 
random variables $\{K_\b(R_{z}(\o))\}_{z\in \G_{y,x}}$ entering in the definition (\ref{w1}) of
$W_\o$. However, thanks to the tree structure of our graph, it is possible
to introduce a modified weight $\wt W_\o(\G_{y,x}):=\prod_{z\in
  \G_{y,x}}\psi_z(\o)$ in such a way that:
  \begin{itemize}
  \item the random variables $\{\psi_z(\o)\}_{z\in \G_{y,x}}$ are
  independent;
\item for each $\o$, $W(\G_{y,x},\o)\leq \wt W(\G_{y,x},\o)$.
  \end{itemize}
We now describe how we construct the modified weights.

To begin with, we fix some notation: $\partial B_\ell$ 
stands for the (deterministic)
set of vertices $x$ such that $d(x)=\ell$. We also use $\G_x$ for the
unique path from the
root to $x$. To simplify the notation we
define $W_\o(x)=0$ if $x\in\partial B_\ell\setminus D_\ell(\o)$.
%, \ie if $x$
%is itself an obstacle. 
Next, for
every $x\in\partial B_\ell$ and for every vertex $z\in\G_x$ we denote by $\D_z$ the set
of all children $y$ of $z$ such that $y\notin \G_x$. Clearly,
$|\D_z|=b-1$.  We say that $z\in\G_x$ is {\em regular} if $R_y\leq \e$
for every $y\in\D_z$. 
\smallno

Let now $u>0$ be a small parameter to be fixed
later and suppose that $z_1,z_2,\dots,z_k$ are consecutive regular
sites on $\G_x$, ordered in such a way that $d(z_j,r)=d(z_{j-1},r)-1$.
Let also $z_0$ denote the children of $z_1$ along $\G_x$.  As in the
proof of Lemma \ref{R_a} we may assume $h\geq -(b-1)+a$ without loss of
generality.  Since $z_1$ is regular, using $|F_\b|_\infty\leq \e^{-1}$,
$F_\b(x)\leq 2x$, for any $x\geq \e$, from (\ref{w3}) we have
$$
R_{z_1}(\o) \leq \e^{a-(b-1)} F_\b(R_{z_0}(\o)) \prod_{y\in\D_{z_1}}F_\b(R_{y}(\o))
\leq \e^{a-b}(2\e)^{b-1} = \e^{a-1}2^{b-1}\,.
$$
Therefore, by proceeding inductively, for any $k$ we have
\begin{align}
R_{z_{k}}(\o) &\leq \e^{a-(b-1)} F_\b(R_{z_{k-1}}(\o)) \prod_{y\in\D_{z_k}}F_\b(R_{y}(\o))
\nonumber\\
&\leq \e^{a-(b-1)}\e^{(k-1)a-1}2^{(k-1)b}(2\e)^{b-1} = \e^{ka-1}2^{kb-1}\,.
\la{rk1}
\end{align}
Suppose now $k\geq k_0:=\inte{\frac 4a}$. The above argument shows that
$R_{z_k}\leq u\e$, provided $\e\leq\e_0(a,b,u)$.  A simple computation
gives $K(\a\e)\leq \a$ for every $\a>0$, so that $K_\b(R_{z_k})\leq u$.
%
%\begin{definition}

We shall say that $z\in\G_x$ is {\em good} if $z$ is regular and the
number of consecutive regular vertices immediately below $z$ along
$\G_x$ is larger or equal to $k_0-1$. Otherwise we say that $z$ is {\em
  bad}.   
%\end{definition}

The estimate (\ref{rk1}) therefore implies that
$K_\b(R_{z}(\o))\leq u$ whenever $z$ is good.  Since $K_\b\leq 1$ and
recalling that $W(x)=0$ if $x$ is not connected to the root in $T(\o)$
we may write 
\begin{equation}
W(x)\leq \wt W(x):=\prod_{z\in\G_x} \psi_z(\o)\,,
\quad\, \psi_z(\o) :=
\begin{cases}
u  & {\o_z=+1\,,\; z \; {\rm is} \;\text{\em good}}\\
1 & {\o_z=+1\,,\; z\;{\rm is} \;\text{\em bad}}\\
0 & {\o_z=-1}
\end{cases}\,
\la{wpsiw}
\end{equation}

%                             %%%%%%%%%%    PICTURE !@!@!!   %%%%%%%
\begin{figure}[h]%t= top, b = bottom, h = here
\centerline{
\psfrag{x}{$x$}
\psfrag{z1}{$z_1$}
\psfrag{z2}{$z_2$}
\psfrag{zk}{$z_k$}
\psfrag{r}{$r$}
%\psfrag{la}{$\La\;\rightarrow$}
%\psfrag{.}{{\Large $\cdot$}} \psfrag{x}{$x$}
\psfig{file=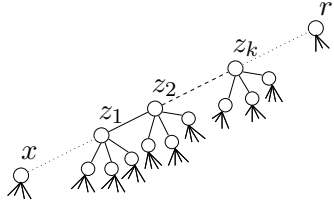,height=1in} }
\caption{$k$ consecutive regular sites on the path $\G_x$
in the case $b=4$.} 
\label{fig4}
\end{figure}                    %%%%%%%%%%    PICTURE !@!@!!   %%%%%%%
We now claim that there exist 
$C_1<\infty$ such that for every $\ell\geq 0$, for all $x\in\partial B_\ell$:
\begin{equation}
\bbE_p \bigl[\wt W(x)\bigr]\leq \,C_1\,(2u)^\ell\,,
%\nep{-m\ell}\,,
\la{ewtw}
\end{equation}
%for every $a>0$, $b\geq 2$ and for any $u>0$ there exist 
%constants $C_1<\infty$, $\b_0<\infty$ and $p_0<1$, 
%such that for any $\b\geq \b_0$, $h\geq -h_c(\b)+a$, $p\geq p_0$ and for every $x\in\partial B_\ell$:
%\begin{equation}
%\bbE_p \wt W(x)\leq \,C_1\,(2u)^\ell\,,
%\nep{-m\ell}\,,
%\la{ewtw}
%\end{equation}
%for all $\ell\geq 0$.
From (\ref{wpsiw}) we see that 
\begin{equation}
\wt W(x)\leq u^{\ell - n_x(\o)}\,.
\la{bad}
\end{equation} 
where $n_x(\o)$
stands for the number of {\em bad} vertices in $\G_x$. 
Define now, for every $z\in\G_x$, $\chi_z(\o)=0$ if $z$ is regular 
and $\chi_z(\o)=1$ otherwise. Note that, by construction, these are i.i.d.\ Bernoulli 
random variables. A simple deterministic bound on $n_x$ is given by 
$$n_x(\o)
\leq k_0 \left(1+\sum_{z\in\G_x}\chi_z(\o)\right)\,.$$ 
From Lemma \ref{R_a} we know that the probability of being irregular, 
for any given $z\in\G_x$, is less than $\d_1:=(b-1)\d + (b-1)(1-p)$. 
Let us choose $\d$ in Lemma  \ref{R_a} and $p<1$ such that 
$\d_1\leq u^{k_0}$. We then have
\begin{equation}
\bbE_p\left[\wt W(x)\right]\leq 
u^{\ell-k_0} \bbE_p\left[u^{- k_0\chi_z}\right]^\ell
 \leq u^{\ell-k_0}(1+\d_1u^{-k_0})^\ell \leq u^{\ell-k_0} 2^\ell\,.
\la{d1u}
\end{equation}
%From the independence of the variables $\chi_z$, $z\in\G_x$, (\ref{wtw}) and (\ref{d1u})
%imply $\bbE_p \wt W(x)\leq 2^\ell u^{\ell-k_0}$. 
The claim (\ref{ewtw}) then follows by taking $C_1=u^{-k_0}$.

\smallskip

We are ready to prove the exponential moment estimate (\ref{expw}). 
For any integer $k$ we define 
\begin{equation}
M_k = \left(\sum_{x\in \partial B_\ell}\wt W(x) \right)^k = 
\sum_{x_1,\dots,x_k\in \partial B_\ell}\wt W(x_1)\cdots \wt W(x_k)\,.
\la{wk1}
\end{equation}
We claim that 
\begin{equation}
\bbE_p \left[ M_{k}\right] \leq  C_2^{k} k!\,,\quad \;k=1,2,\dots
\la{expw1}
\end{equation}
for some constant $C_2<\infty$.
Note that the result (\ref{expw}) is an immediate consequence of (\ref{expw1}) since
the l.h.s.\ in (\ref{expw}) is bounded by 
$$
\bbE_p\left[\exp{(tM_1)}\right] = 
1 + \sum_{k=1}^\infty \frac{t^{k}}{k!}\bbE_p M_{k}
\leq \frac{1}{1-C_2t}\,.
$$
Let $x_1,\dots,x_k\in \partial B_\ell$ be given as in a generic term in the
sum in (\ref{wk1}). These points may be ordered by the 
lexicographic rule to obtain 
the ordered set $\tilde x_1\leq \tilde x_2\leq \dots \leq \tilde x_k$. 
Call 
$\tilde x_0$ and $\tilde x_{k+1}$ the absolute leftmost and, respectively, the absolute rightmost
vertex in $\partial B_\ell$. 
%and 
%$\tilde x_{k+1}$ the last vertex in $D_\ell(\o)$. 
Below we use $[\tilde x_{j-1},\tilde x_j)$ to denote the set of
vertices $y\in\partial B_\ell$ such that $y$ 
is larger or equal to $\tilde x_{j-1}$ 
but strictly less than $\tilde x_j$, with the agreement that, when $j=k+1$,
the set 
$[\tilde x_{k},\tilde x_{k+1})$ also includes the end point $\tilde x_{k+1}$. 
With these notations we can write
\begin{equation}
M_{k+1}(\o) = 
%\sum_{x_1,\dots,x_k\in D_\ell(\o)}W(x_1)\cdots W(x_k)
%\sum_{y\in D_\ell(\o)}W(y) \nonumber\\
%& =
\sum_{j=1}^{k+1}\,
\sum_{x_1,\dots,x_k\in \partial B_\ell}\wt W(x_1)\cdots \wt W(x_k)
\sum_{y\in [\tilde x_{j-1},\tilde x_j)} \wt W(y) 
\la{wk2}
\end{equation}

%                             %%%%%%%%%%    PICTURE !@!@!!   %%%%%%%
\begin{figure}[h]%t= top, b = bottom, h = here
\centerline{
\psfrag{x0}{$\tilde x_0$}
\psfrag{x1}{$\tilde x_1$}
\psfrag{x2}{$\tilde x_2$}
\psfrag{x3}{$\tilde x_3$}
\psfrag{x4}{$\tilde x_4$}
\psfrag{x5}{$\tilde x_5$}
\psfrag{r}{$r$}
\psfrag{y}{$y$}\psfrag{zy}{$z_y$}
%\psfrag{la}{$\La\;\rightarrow$}
%\psfrag{.}{{\Large $\cdot$}} \psfrag{x}{$x$}
\psfig{file=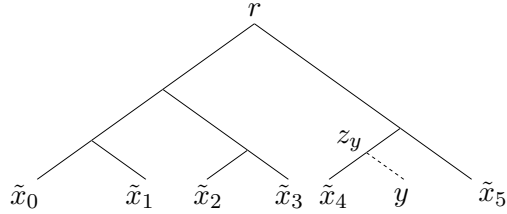,height=1in,width=2.5in} }
\caption{Schematic picture of the tree $T_k$ in the case $k=4$. 
Here $y\in[\tilde x_4,\tilde x_5)$.} 
\label{fig5}
\end{figure}                    %%%%%%%%%%    PICTURE !@!@!!   %%%%%%%
Consider now 
a given $y\in [\tilde x_{j-1},\tilde x_j)$. 
Let $T_k$ denote the subtree determined by the
union of all paths $\G_{r,x_i}$, $i=1,\dots,k$.
Let $d(y,T_k)$
denote the distance from $y$ to $T_k$
and write $z_y$ for the closest ancestor of $y$ on $T_k$
(characterized by $d(z_y,y)=d(y,T_k)$). Clearly, we can estimate
$\wt W(y)\leq \wt W(\G_{z_y,y})$, where $\wt W(\G_{z_y,y})=\prod _{z\in\G_{z_y,y}}\psi_z$.
Now, by construction, the random variable $\wt W(\G_{z_y,y})$ is independent of all the 
weights $\wt W(x_i)$ except for the variables $\psi_z$ where $z$ is 
either $z_y$ or one of the 
$k_0-1$ consecutive vertices just above $z_y$. Let us call $A_y$ this set of vertices. 
Restricting to the event that 
$x_1,\dots,x_k\in D_\ell(\o)$ we 
can estimate 
$$
%\wt W(x_1) \cdots \wt W(x_k)\leq \wt W(x_1) \cdots \wt W(x_k)/
%\left(
\prod_{z\in A_y} \left(\psi_z(\o)\right)^{-1}
\leq u^{-k_0}\,.  %\wt W(x_1) \cdots \wt W(x_k)\,.
$$ 
Therefore we have 
\begin{equation}
\bbE_p \left[M_{k+1}\right] \leq  \,u^{-k_0}\sum_{x_1,\dots,x_k\in\partial B_\ell}\,
\bbE_p\left[W(x_1)\cdots W(x_k)\right]
\sum_{j=1}^{k+1}\sum_{y\in [\tilde x_{j-1},\tilde x_j)} 
\bbE_p \left[\wt W(\G_{z_y,y})\right]
\la{w5}
\end{equation}
Clearly, for every integer $d$ there at most 
$b^d$ vertices $y$ such that $d(y,T_k)=d$. Therefore,
choosing $u\leq 1/(4b)$, 
from (\ref{ewtw}), %there exists $C<\infty$ such that 
for every pair $\tilde x_{j-1},\tilde x_j$, we have
\begin{equation}
\sum_{y\in [\tilde x_{j-1},\tilde x_j)}  \bbE_p \left[\wt W(\G_{z_y,y}) \right]
\leq 2C_1\,.
\la{ub}
\end{equation}
From (\ref{w5}) and (\ref{ub}), setting $C_2:=2u^{-k_0}C_1$ we obtain
\begin{equation}
\bbE_p \left[M_{k+1}\right]\leq  C_2(k+1)\,\bbE_p \left[M_{k}\right]  \,,
\la{w6}
\end{equation}
so that, for every $k\geq 1$ we can estimate 
$\bbE_p \left[M_{k}\right] \leq  C_2^{k}k!$ as claimed in (\ref{expw1}). 
\qed

\bigskip 

We turn to the setting of statement \emph{b*}. \begin{lemma}
\la{weights_b}
For any $p>\frac12$, there exist $\b_0<\infty$ and $b_0\in\bbN$, 
such that (\ref{expw}) holds for any $\b\geq \b_0$, $b\geq b_0$, $h=0$.
\end{lemma}
\proof
The proof is essentially the same as that of Lemma \ref{weights_a},
but we have to modify the definition of good and bad vertices. 
Given $x\in\partial B_\ell$ and $z\in\G_x$ we write as before $\D_z$
for the set of children of $z$ lying outside of the path $\G_x$. 
We write also $m(\o)$ for the number of $y\in\D_z$ such that $\o_y=-1$.
As in the proof of Lemma \ref{R_b}, we may use (\ref{expb})
to estimate this quantity: there exist $a_1,a_2>0$ such that 
$\bbP_p\left( m \geq (1-2a_1)b/2 \right)
\leq \nep{-a_2 b}$, for all sufficiently large $b$.

\smallskip

Here the definition of good vertices goes as follows.
We say that $z$ is \emph{good}
if $m\leq (1-2a_1)b/2$ and if all the vertices $y\in\D_z$ such
that $\o_y=+1$ satisfy $R_y\leq \e$. Clearly, if $z$ is good, from (\ref{w3}) 
we must have
$$
R_z\leq \e^{-m}(2\e)^{b-1-m}\e^{-1} \leq 2^{b-1} \e^{2(a_1b-1)}\,.
$$
In particular, for any $u>0$ we find $b_0$ and $\b_0$ such that
for all $b\geq b_0$ and $\b\geq \b_0$ we have $R_z\leq u\e$.  
As in the proof of Lemma \ref{weights_a} we therefore have that $K_\b(R_z)\leq u$ 
whenever $z$ is good. Now we can define individual weights $\psi_z$ exactly as in (\ref{wpsiw})
and, as before, we can estimate $W(x)\leq \wt W(x)$.
To establish the analog of (\ref{ewtw}) we simply observe that 
$\psi_z$ are i.i.d.\ random variables 
with the present definition of good vertices. 
Moreover, from Lemma \ref{R_b} we easily infer that
\begin{equation}
\bbP_p\left[\,z \; \text{is bad}\right] \leq \nep{-c\,b}
\la{nepb}
\end{equation}
for some $c>0$. Using this we have, see (\ref{d1u}) 
\begin{equation}
\bbE_p\left[\wt W(x)\right]\leq (2u)^\ell\,,
%u^{\ell}(1+ \nep{-c\,b} u^{-1})^\ell \leq u^{\ell} 2^\ell\,.
\la{dbu}
\end{equation}
as soon as $\nep{-c\,b} u^{-1}\leq 1$. The rest of the proof goes now exactly as
in Lemma \ref{weights_a}. The estimate (\ref{w5}) is 
actually simplified by the fact that we only need to remove
the vertex $z_y$, so that the factor $u^{-k_0}$ is now replaced by $u^{-1}$. 
In particular, (\ref{w6}) now holds with the constant $C_2 = 2u^{-1}$. 
\qed
\bigskip 

It remains to prove (\ref{expw}) in the setting of statement \emph{c*}. 
\begin{lemma}
\la{weights_c}
For any $p>0$, there exist $\b_0<\infty$ and $b_0\in\bbN$, 
such that (\ref{expw}) holds for any $\b\geq \b_0$, $b\geq b_0$, $h=+h_c(\b)$.
\end{lemma}
\proof
As in the proof of Lemma \ref{weights_c} we call $m(\o)$ the number of obstacles 
among the children in $\D_z$ for a given vertex $z$.
We shall use the analog of estimate (\ref{expc}). Letting $a_1$ and $a_2$ be the 
parameters appearing there, the vertex $z$ 
is now declared \emph{good} if $m(\o)$ satisfies $m(\o)\leq (1-a_1)b$ and 
all the vertices $y\in\D_z$ such that $\o_y=+1$ satisfy $R_y\leq \e$. 

\smallskip

With this definition of good vertices, 
the bounds of Lemma \ref{R_c} now show that 
$$
\bbP_p\left[\,z \; \text{is bad}\right] \leq \nep{-c\,b} 
$$
for some $c>0$ and all sufficiently large $b$.
On the other hand, reasoning as in the proof of Lemma \ref{R_c}, 
an application of (\ref{w3}) gives that if $z$ is good then 
$$
R_z\leq \e^{b-2}\e^{-m}(2\e)^{b-1-m}\e^{-1}\leq 2^{b-1}
\e^{2a_1b -4}\,.
$$
Given $u>0$ we then find $b_0,\b_0$ such that $K_\b(R_z)\leq u$ 
whenever $z$ is good, as soon as $b\geq b_0$ and $\b\geq \b_0$. 
The rest goes exactly as in the proof of Lemma \ref{weights_b}
above. 
\qed

\subsection{Poincar\'e and 
Logarithmic Sobolev inequalities among obstacles}
Recall the definition (\ref{eq:cgapcsob}) of the constants $\cgap$ and $\csob$.
We shall focus here on the case of the measure $\mu=\mu_{L,\o}^+$, \ie
the Gibbs measure with plus boundary condition below a certain level $L$ among the
obstacle environment $\o$. We shall write $\cgap(L,\o)$ and $\csob(L,\o)$ for the associated constants.
Here $\o$ will be distributed according to Bernoulli(p) measure.

\smallskip
 
It is well known that, in general, $c_{\rm gap} \geq 2 c_{\rm sob}$. On the other hand,
for trees, a useful inequality established in \cite{MaSiWe} states
that $c_{\rm gap} \leq O(\log n) c_{\rm sob}$, where $n$ is the cardinality of 
the tree. In particular, Theorem 5.7 in \cite{MaSiWe} in our setting implies
that for every $b$ and every $\beta$ 
there exists a constant $C<+\infty$ such that for every $\o\in\O$ and for 
every $L$  
\begin{equation}
c_{\rm gap}(L,\o)\leq C\,L\, c_{\rm sob}(L,\o)
\la{crude}
\end{equation}

\smallskip

Our main result here is an almost sure polynomial bound on $c_{\rm sob}(L,\o)$:
There exists a constant $\zeta<\infty$ such that 
\begin{equation}
c_{\rm sob}(L,\o)\geq L^{-\zeta} 
\la{polybound}
\end{equation}
holds $\bbP_p$--a.s.\ for $L$ sufficiently large. Here 
the constant $\zeta$ will depend on the parameters 
$a,b,\b$ in case {\em a*}, on $b,\b$ and $p$ in cases \emph{b*,c*}.
\begin{theorem}
\label{theosob}
\medno 
\begin{enumerate}[a*)]
\item For every $a>0$, $b\geq 2$,
there exists $p<1$ and $\b_0$ such that
(\ref{polybound}) holds 
for all $\b\geq \b_0$ and $h\geq -h_c(\b,b) + a$.\medno \item
For every $p>\frac12$, there exist $b_0\in\bbN$ and $\b_0<\infty$ such that 
(\ref{polybound}) holds 
for $h=0$, $b\geq b_0$, $\b\geq \b_0$.\medno \item
For every $p>0$, there exist $b_0\in\bbN$ and $\b_0<\infty$ such that 
(\ref{polybound}) holds 
for $b\geq b_0$, $\b\geq \b_0$ and $h=+h_c(\b)$.\end{enumerate}
\end{theorem}
\proof
We will carry out the proof of the three statements simultaneously.
Indeed, the key estimate we need is the exponential integrability (\ref{expw}),
which holds in all cases under consideration as worked out in 
Lemma \ref{weights_a}, Lemma \ref{weights_b} and Lemma \ref{weights_c}.

\smallskip

Thanks to the deterministic bound (\ref{crude}) it suffices to
prove the claim (\ref{polybound}) 
with $c_{\rm sob}(L,\o)$ replaced by $c_{\rm gap}(L,\o)$.
We fix a length scale $\ell_1$ much smaller than $L$.
For each vertex $x\in T_L(\o)$, let $B_{x,\ell_1}\subset T_L(\o)$
denote the subtree (or ``block'') of depth $\ell_1 - 1$ rooted at $x$.
In this way $B_{x,\ell_1}$ consists of $\ell_1$ levels and we understand that 
if $x$ is $k<\ell_1$ levels from the bottom of $T_L(\o)$ then 
$B_{x,\ell_1}$ has only $k$ levels. In the end we will choose $\ell_1=C\log L$
for some sufficiently large constant $C$.
% depending on the various parameters. 
We define
the Dirichlet form of the so-called ``block--dynamics''
$$
\cD_{\ell_1,L,\o}(f) = \sum_{x\in T_L(\o)} \mu^+_{L,\o}
\left[\var_{B_{x,\ell_1}}(f)\right]\,.
$$
A standard argument relating the spectral gap 
of the heat--bath dynamics to the spectral gap of the block--dynamics
(see e.g.\ \cite{Mar}) shows that, since there are at most $\ell_1$ blocks
containing a given vertex $x$, we have  
\begin{equation}
c_{\rm gap}(L,\o)\geq \frac1{\ell_1}\,\left\{
\min_{\tau,x}c_{\rm gap}(\mu_{B_{x,\ell_1}}^\tau)\right\}\,
\inf_{f} \frac{\cD_{\ell_1,L,\o}(f)}{\var_{\mu^+_{L,\o}}(f)}\,,
\la{gap1}
\end{equation}
where $c_{\rm gap}(\mu_{B_{x,\ell_1}}^\tau)$ denotes the spectral gap of the 
heat bath dynamics on the block $B_{x,\ell_1}$ with boundary condition $\tau$
(and $\tau$ is assumed to be compatible with the obstacle realization, i.e.\
$\tau\in\cB_\o$):
$$
c_{\rm gap}(\mu_{B_{x,\ell_1}}^\tau)=\inf_{f}\,
\frac{\cD_{\mu_{B_{x,\ell_1}}^\tau}(f)}{\var_{\mu_{B_{x,\ell_1}}^\tau}(f)}
\,.
$$ 
In general trees, according to Theorem 1.4 in \cite{BKMP}, one has a 
lower bound on $c_{\rm gap}$ of order $n^{-\zeta}$ uniformly over 
the boundary condition, where $n$ is the cardinality of the tree
and $\zeta<\infty$ is a constant depending on the parameters $b,\b,h$.
In particular this implies that 
for all $\o\in\O$ and for all sufficiently large $\ell_1$  
\begin{equation}
\frac1\ell_1\,\min_{\tau,x}c_{\rm gap}(\mu_{B_{x,\ell_1}}^\tau)\geq 
b^{-2\zeta\ell_1}\,.
\la{gap2}
\end{equation}
Let $c_{\rm gap}(\ell_1,L,\o)$ denote the spectral
gap of the block--dynamics, \ie the infimum appearing in (\ref{gap1}).
So far we have obtained the deterministic bound
\begin{equation}
c_{\rm gap}(L,\o)\geq b^{-2\zeta\ell_1}\,c_{\rm gap}(\ell_1,L,\o)\,.
\la{gap3}
\end{equation}
Next we make a deterministic estimate on $c_{\rm gap}(\ell_1,L,\o)$.
To this end we use the method of \cite{MaSiWe},
combined with the results we obtained in previous subsections. 
Given $r\in(0,1)$, we say that $\mu^+_{L,\o}$ is $(\ell_1,r^{\ell_1})$--mixing if 
%there exists $a>0$ such that 
for every $x\in T_L(\o)$
\begin{equation}
\var_{\mu^+_{L,\o}}\left(\mu^+_{L,\o}(\si_x\tc \si_{D_{x,\ell_1}(\o)})\right)
\leq r^{\ell_1} \var_{\mu^+_{L,\o}}(\si_x)\,,
\la{gap4}
\end{equation}
where $\mu^+_{L,\o}(\si_x\tc \si_{D_{x,\ell_1}(\o)})$ denotes the
conditional expectation of $\si_x$ given the values of $\si$
on $D_{x,\ell_1}(\o)$, the set of descendants of $x$ at distance $\ell_1$.
A simple computation shows that (\ref{gap4}) is actually equivalent to 
the variance mixing condition ${\rm VM}(\ell_1,\eps)$, with $\eps=r^{\ell_1}$,
introduced in \cite{MaSiWe}. In particular, 
Theorem 3.2 in \cite{MaSiWe} implies 
that 
\begin{equation}
c_{\rm gap}(\ell_1,L,\o)\geq \frac 14\,,
\la{gap40}
\end{equation} 
for $\ell_1\geq \ell_0$,
for some finite $\ell_0=\ell_0(r)$ as soon as $\mu^+_{L,\o}$ is $(\ell_1,r^{\ell_1})$--mixing.
The conclusion of the theorem therefore follows from (\ref{gap3}) and (\ref{gap40})
if we can prove that $\mu^+_{L,\o}$ is $(\ell_1,r^{\ell_1})$--mixing $\bbP_p$-a.s.\
for some $r<1$, when $\ell_1=C\log L$, 
with some $C<\infty$, for all 
sufficiently large $L$.
To prove this we observe that, setting $g(\si_x):=
 \mu^+_{L,\o}\left[\mu^+_{L,\o}(\si_x\tc \si_{D_{x,\ell_1}(\o)})\tc \si_x\right]$,
we may write
\begin{align}
&\var_{\mu^+_{L,\o}}
\left(\mu^+_{L,\o}(\si_x\tc \si_{D_{x,\ell_1}(\o)})\right)
= \Cov_{\mu^+_{L,\o}}\left(\si_x,g(\si_x)
\right)\nonumber\\
&= 2\, \mu^+_{L,\o}(\si_x=+1)\,
\mu^+_{L,\o}(\si_x=-1)\, \left[g(+1)-g(-1)
\right] \leq \frac12 \left[g(+1)-g(-1)
\right]\,.
\la{gap01}
\end{align}
Let $\nu$ denote a coupling of the measures $\mu^+_{L,\o}(\cdot\tc \si_x=+1)$
and $\mu^+_{L,\o}(\cdot\tc \si_x=-1)$. We then write
\begin{align*}
& g(+1)-g(-1)  = 
\sum_{\tau,\eta}\nu(\tau,\eta)\left[
\mu^+_{L,\o}(\si_x \tc \tau_{D_{x,\ell_1}(\o)})
- \mu^+_{L,\o}(\si_x\tc \eta_{D_{x,\ell_1}(\o)})\right]\\
&\quad \leq  \sum_{\tau,\eta}\nu(\tau,\eta)
\sum_{y\in D_{x,\ell_1}(\o)} 1_{\tau_y\neq\eta_y}\left[
\mu^+_{L,\o}(\si_x \tc ([\tau\eta]^{y,+})_{D_{x,\ell_1}(\o)})
- \mu^+_{L,\o}(\si_x\tc ([\tau\eta]^{y,-})_{D_{x,\ell_1}(\o)})\right]
\end{align*}
with $[\tau\eta]^{y,\pm}$ denoting the interpolation 
between $\tau$ and $\eta$, i.e.\ the configuration such that, 
using lexicographic order on $D_{x,\ell_1}(\o)$, 
$([\tau\eta]^{y,\pm})_z=\tau_z$, 
$z<y$, and $([\tau\eta]^{y,\pm})_z=\eta_z$, $z>y$,
while $([\tau\eta]^{y,\pm})_y=\pm 1$.
Recall now the definition (\ref{k1}) of the function $K_\b$
and set 
$$
\g:=\sup_{x>0}K_\b(x) = \tanh \b \,.
$$
Since $[\tau\eta]^{y,+}$ and $[\tau\eta]^{y,-}$
differ only at $y\in D_{x,\ell_1}(\o)$, reasoning as in (\ref{totvar1})
and (\ref{totvar2}) we estimate
\begin{equation}
\mu^+_{L,\o}(\si_x \tc ([\tau\eta]^{y,+})_{D_{x,\ell_1}(\o)})
- \mu^+_{L,\o}(\si_x\tc ([\tau\eta]^{y,-})_{D_{x,\ell_1}(\o)})
\leq 2\,\g^{\ell_1}
\la{gap6}
\end{equation}
From (\ref{ham}) we then obtain
\begin{equation}
g(+1)-g(-1) \leq 2\,\g^{\ell_1}\,\nu\left(|\tau-\eta|_{x,\ell_1}\right)
\leq 2\,\g^{\ell_1} \sum_{y\in D_{x,\ell_1}(\o)} W(\G_{x,y})
\la{gap7}
\end{equation}
Since $\d=\d(\b,h,b):=\var_{\mu^+_{L,\o}}(\si_x)>0$, 
from (\ref{gap4}), (\ref{gap01}) and (\ref{gap7})
we see that
\begin{align}
&\bbP_p\left( (\ell_1,r^{\ell_1}) {\rm - mixing \;does \;not \;hold}
\right)
\leq 
\bbP_p\Big(\exists B_{x,\ell_1}: \sum_{y\in D_{x,\ell_1}} W(\G_{x,y})
\geq \d\,\g^{-\ell_1}\,r^{\ell_1} \Big)
\nonumber\\
& 
\hskip5cm\leq \sum_{x:\; d(x)\leq L}
\, 
\bbP_p\Big(\sum_{y\in D_{x,\ell_1}} W(\G_{x,y})
\geq \d\,\g^{-\ell_1}\,r^{\ell_1} \Big)\,.
\la{gap8}
\end{align}
For every $x$ we can use the bound (\ref{expw}), so that 
(\ref{gap8}) yields 
\begin{equation}
\bbP_p\left((\ell_1,r^{\ell_1}){\rm - mixing \;does \;not \;hold}
\right)\leq 2\, b^L\, \exp{(-t_0\d \,r^{\ell_1}\g^{-\ell_1})}\,.
\la{gap9}
\end{equation}
Setting e.g.\ $r = \sqrt \g$, $\ell_1=C\log L$ with $C$ sufficiently large
we see that by the Borel Cantelli lemma we have $(\ell_1,r^{\ell_1})$--mixing
$\bbP_p$-a.s.\ for all $\ell_1$ large enough.
This concludes the proof of the theorem. \qed

\medno

{\bf Remark}. 
One may wonder whether the result of Theorem \ref{theosob} captures
  the true behavior of the logarithmic Sobolev constant in presence of
  the random realization of obstacles or whether, instead, it only provides a
  pessimistic bound. As we show below, as soon as $\b$ is large enough
  (actually larger than the spin--glass critical point for the pure
  Ising model on $\Tree^b$ \cite{BKMP}), in all the three cases described in
  the theorem, there exists a set $\O_0$ of obstacles realizations of uniformly positive 
probability, such that for every $\o\in\O_0$ 
  the spectral gap and a fortiori the logarithmic Sobolev constant must
  shrink to zero at least as fast as $L ^{-\zeta'}$ for some deterministic
  exponent $\zeta'>0$. A quick sketch of the proof of this fact for the
  setting $(a^*)$ and e.g. $b=2$ goes as
  follows.

  Pick a vertex $x\in T(\o)$ with $d(x)\leq \frac L2$ and denote by
  $T_{x,\ell}$ the finite sub-tree of $\Tree^b$ rooted at $x$ with
  $\ell=\d(M+1)\log L$ levels, $\d\ll 1,\, M\gg 1$ but $\d M \ll 1$.
  Assume that $T_{x,\ell}$ is free of obstacles, group together the
  sites of $\partial T_{x,\ell}$ into equal blocks according to their
  common ancestor in $T_{x,\ell}$ at level $\d\log L +1$ and order the
  blocks from left to right. Then impose that all vertices inside the
  odd blocks are obstacles while all sites $z$ inside even blocks are
  not obstacles and the corresponding $R_z(\o)$ satisfies $R_z(\o)\le
  \e$ (as usual $\e=\nep{-2\b}$).  Because of Lemma
  \ref{R_a} the obstacles realizations that obey the above
  specifications have probability larger than $\nep{-c|T_{x,\ell}\cup
    \partial T_{x,\ell}|} = \nep{-c L^\a}$ for a suitable constant
  $c=c(p,\b)$ and $\a=\d(M+1)\log 2$. 
Since $\d M \ll 1$ the probability of finding a vertex $x$ with the
  above properties converges to one as $L\to \infty$. 

Consider now the common ancestor
  $y$ at level $\d\log L$ of two odd-even neighboring blocks. It follows
  immediately from the recursion (\ref{w3}) and the assumptions we made
  on the even/odd blocks, that $|R_y-1|\leq \nep{-c\d M\log L}$ for a
  suitable constant $c$. In turn, if $M$ is large enough, that implies
  that the marginal of the Gibbs measure $\m_{L,\o}$ on the finite
  sub-tree rooted at $x$ with now $\d\log L$ levels, has a bounded
  (independently of $L$) relative density with respect to the Ising
  Gibbs measure on the same tree with \emph{free boundary conditions} on
  its leaves. Since the latter has a spectral gap (and a fortiori a
  logarithmic Sobolev constant) smaller than $b^{-a(\b)\d\log L}$ for
  some positive $a(\b)$, we conclude that for the obstacles realizations
  satisfying the previous conditions, $\csob(\mu^+_{L,\o})\leq L^{-\zeta'}$
  for some $\zeta'=\zeta'(\b,\d,b)$.  

%\end{remark}
%}

\section{Proof of Claims (1)--(4)}\label{claims}
The results of the previous section allow us to fill the gaps
in the proof of Theorem \ref{theotree}. We refer to section \ref{main}
for the setting and the notation. 

\medno
{\bf Claim 3}. 
Recall the definition of the random trees $T_{z,2\ell}(\o)$ and the associated measure 
$\nu^{+,+}_{z,2\ell}$, where $z$ is a vertex at level $\ell$. 
We have to estimate $\min_z\cgap(\nu_{z,2\ell}^{+,+})$. Since the 
b.c.\ above $z$ can only affect this quantity by a 
constant factor (depending on $\b$), we may replace $\nu_{z,2\ell}^{+,+}$
by the measure $\nu_{z,2\ell}^+$ with free b.c.\ above $z$.
At this point, for each $z$ we are exactly 
in the setting of Theorem \ref{theosob}, with $L=2\ell$. 
As we have seen in the proof of that theorem (see (\ref{gap9})),
we can prove that
there exists $\zeta<\infty$ such that 
\begin{equation}
\bbP_p\left(\csob(\nu_{z,2\ell}^+)\leq \ell^{-\zeta}\right)\leq \nep{-\ell^2}\,,
\la{c31}
\end{equation}
for all sufficiently large values of $\ell$.
Since the number of $z$ such that $d(z)=\ell$ is $b^\ell$,
the claim follows from the Borel--Cantelli lemma. \qed

\medno
{\bf Claim 4}. Observe that it is sufficient to prove
\begin{equation}
\bbP_p\left(\nu_{r_x,2\ell}^{+,+}(\si_x = +1) -
\nu_{r_x,\infty}^{-,+}(\si_x = +1)\geq (3b)^{-3\ell}
\right)\leq \nep{-\ell}\,.
\la{c41}
\end{equation}
Recall that $\nu_{r_x,2\ell}^{-,+}$ stands for the Gibbs measure on $T_{r_x,2\ell}(\o)$ 
with $-$ b.c.\ above $D_\ell(\o)$ and $+$ b.c.\ below $D_{2\ell}(\o)$. 
Observe that, by (\ref{totvar2}) we have
\begin{equation}
\nu_{r_x,2\ell}^{+,+}(\si_x = +1) -
\nu_{r_x,2\ell}^{-,+}(\si_x = +1) = W(\G_{r_x,x})\,.
\la{c42}
\end{equation}
As in the proof of (\ref{expw}) the 
expectation of the latter expression is estimated by $\bbE_p\wt W(\G_{r_x,x})\leq (2u)^{\ell/2}$,
since $x$ satisfies $d(x,r_x)\geq \ell/2$ (see also Lemma \ref{weights_a},
Lemma \ref{weights_b} and Lemma \ref{weights_c}). 
Then by Markov's inequality
\begin{equation}
\bbP_p\left(\nu_{r_x,2\ell}^{+,+}(\si_x = +1) -
\nu_{r_x,2\ell}^{-,+}(\si_x = +1)\geq \frac12(3b)^{-3\ell}
\right)\leq \frac12\,\nep{-\ell}
\la{c420}
\end{equation}
provided $u$ is small enough. We turn to an estimate of the difference
$\nu_{r_x,2\ell}^{-,+}(\si_x = +1) - \nu_{r_x,\infty}^{-,+}(\si_x =
+1)$. Let $E_\ell(\o)$ denote the (Ising--model) event that there exists
a path $\G$ in $T(\o)$ joining the sets $D_{\frac52 \ell}(\o)$ and
$D_{3\ell}(\o)$ such that $\si_z=-1$ for each $z\in\G$. If the sets
$D_{\frac52 \ell}(\o)$ and $D_{3\ell}(\o)$ are not connected in $T(\o)$
we simply set $E_\ell = \emptyset$. Observe that by monotonicity, for
every $\o$ we have
\begin{equation}
\nu_{r_x,\infty}^{-,+}(\si_x = +1\tc 
E_{\ell}^c ) \geq \nu_{r_x,2\ell}^{-,+}(\si_x = +1)
\la{c43}
\end{equation}
The reason for the above domination is that if there is no path connecting
$D_{\frac52 \ell}(\o)$ and $D_{3\ell}(\o)$ covered by $-$ spins, then there must exist
a cut--set of $T(\o)$, fully contained between level $\frac52 \ell$ and level $3\ell$
covered by $+$ spins, and conditioned on this event $\nu_{r_x,\infty}^{-,+} $
dominates $\nu_{r_x,2\ell}^{-,+}$. 
We then have 
\begin{align}
\nu_{r_x,2\ell}^{-,+}&(\si_x = +1) -
\nu_{r_x,\infty}^{-,+}(\si_x = +1)
\nonumber
\\
&  \leq 
\nu_{r_x,2\ell}^{-,+}(\si_x = +1) -
\nu_{r_x,\infty}^{-,+}(E^c_{\ell})
\nu_{r_x,\infty}^{-,+}(\si_x=+1 \tc E^c_{\ell}) \nonumber\\
&\leq 
\nu_{r_x,\infty}^{-,+}(E_{\ell})
\la{c44}
\end{align}
At the price of a $\b$--dependent factor we may replace $\nu_{r_x,\infty}^{-,+}$
with the measure $\nu_{r_x,\infty}^{+}$ with free b.c.\ above the vertex $r_x$. 
Now we are in the familiar setting of the previous subsections.
To estimate $\nu_{r_x,\infty}^{+}(E_{\ell})$, 
suppose $\o\in\O$ is such that $T(\o)$ contains a given path 
$\G=\{x_0,x_1,\dots,x_h\}$ with $d(x_{j+1})=d(x_j)+1$, $d(x_0)=\frac52 \ell$
and $d(x_h)=3\ell$, $h=\ell/2$. Write $\{\G=-\}$ for the event 
$\{\si_{x_0} = \si_{x_1} = \cdots = \si_{x_h} =-1\}$ and let $q_j(\o)$
denote the probability $\nu_{r_x,\infty}^{+}(\si_{x_{j}}=-1\tc
\si_{x_{j-1}}=-1)$. Clearly we have
\begin{equation}
\nu_{r_x,\infty}^{+}
(\G=-)=
\nu_{r_x,\infty}^{+}(\si_{x_0}=-1)\prod_{j=1}^{h}q_j(\o)\leq \prod_{j=1}^{h}q_j(\o)\,.
\la{ekl3}
\end{equation}
Observe that 
\begin{equation}
q_j(\o)= \e^{-1} R_{x_j}(\o)/(1+\e^{-1} R_{x_j}(\o))\,.
\la{c45}
\end{equation} 
As in the proof of Lemma \ref{weights_a},
Lemma \ref{weights_b} and Lemma \ref{weights_c} we discriminate the vertices along $\G$ in
good and bad vertices. We know that if $z$ is good, then $R_z(\o)\leq u\e$ and therefore,
by (\ref{c45}) we have $q_j(\o)\leq u$. As in (\ref{bad}) we may then estimate
\begin{equation}
\nu_{r_x,\infty}^{+}
(\G=-)\leq u^{\frac\ell 2 -n(\o)}
 \la{c46}
\end{equation}
where $n(\o)$ stands for the number of bad vertices along $\G$. Summing over all possible
paths and estimating as in (\ref{d1u}) and (\ref{dbu}) we arrive at
\begin{equation}
\bbE_p\left[\nu_{r_x,\infty}^{+}
(E_{\ell})\right]\leq 
b^{3\ell}\, (2u)^{\frac\ell 2} \leq 
\frac14 (3b)^{-3\ell}\nep{-\ell}\,, 
\la{ekl5}
\end{equation}
%for all sufficiently large $\ell$,
provided $u$ is suitably small.
From (\ref{c44}), Markov's inequality yields
\begin{equation}
\bbP_p\left(\nu_{r_x,2\ell}^{-,+}(\si_x = +1) -
\nu_{r_x,\infty}^{-,+}(\si_x = +1)\geq 
\frac12(3b)^{-3\ell}
\right)\leq \frac12\,\nep{-\ell}\,.
\la{c47}
\end{equation}
This, together with (\ref{c420}), ends the proof of Claim 4.
\qed

\smallskip

Let us now turn to Claim 1 and Claim 2. 
Here the environment $\o$ is given by a Bernoulli(p) configuration
$\eta$ below level $\ell$ and is deterministically free of obstacles up
to and including level $\ell$, as prescribed by (\ref{rule}).

\medno
{\bf Claim 1}. Let $E_\ell$ be the event that there exists 
a vertex $x$ with $d(x)=\ell$, such that $\si_z=-1$ for every $z\in\G_x$,
\ie if the root is connected to level $\ell$ by a path covered with $-$ spins.
As in (\ref{c43}) and (\ref{c44}) we have
\begin{equation}
0\leq \mu^+(\si_r) - \mu_\o^+(\si_r) \leq \mu^+_\o(E_\ell)\,.
\la{c11}
\end{equation}
Following (\ref{ekl3}) and (\ref{c45}) we estimate
\begin{equation}
\mu^+_\o(E_\ell)\leq \sum_{x:\;d(x)=\ell}\,
\prod_{z\in\G_x} (\e^{-1}R_z(\o))\,.
\la{c12}
\end{equation}
By monotonicity we have $R_z(\o)\leq R_z(\wt \o)$, where $\wt \o$
coincides with $\o$ (and therefore with $\eta$) below level $\ell$ and
is given by a new (independent) Bernoulli(p) configuration $\eta'$ up to
and including level $\ell$. We denote by $\wt\bbE_p$ the expectation
over the random environment $\wt \o$. Here we can apply the machinery
developed in Lemma \ref{weights_a}, Lemma \ref{weights_b} and Lemma
\ref{weights_c}.  Namely, for a suitably small parameter $u>0$, we can
write $\e^{-1}R_z(\wt\o)\leq u$ for every good vertex $z$. Estimating as
in (\ref{c46}) and (\ref{c47}) above we have
\begin{equation}
\bbE_p\left[\mu^+_\o(E_\ell)
\right]\leq 
 \sum_{x:\;d(x)=\ell} \wt \bbE_p\left[ \prod_{z\in\G_x} (\e^{-1}R_z(\wt \o))
\right]\leq b^\ell\,(2u)^\ell \leq \nep{-3\ell}
\,.
\la{c120}
\end{equation}
Therefore
\begin{equation}
\bbP_p\left(\mu^+_\o(E_\ell) > \nep{-2\ell}
\right)\leq \nep{2\ell}\,\bbE_p\left[\mu^+_\o(E_\ell)\right]\leq \nep{-\ell}\,.
\la{c13}
\end{equation}
Thanks to (\ref{c11}) and the Borel--Cantelli lemma, this implies 
the desired estimate. \qed

\medno
{\bf Claim 2}. This is the same as the statement (\ref{polybound})
appearing in Theorem \ref{theosob}, with the difference that now the
environment is deterministically free of obstacles up to and including
level $\ell$. We can therefore repeat the argument used in the proof of
Theorem \ref{theosob} and see that what has to be established here is a
version of the exponential integrability (\ref{expw}) for our new
environment $\o$.  The latter, in turn, relies on the bounds of Lemma
\ref{R_a}, Lemma \ref{R_b} and Lemma \ref{R_c}. Since the ratios $R$ are
monotonic functions of the environment, these estimates can only improve
in the setting considered here and the proof of Claim 2 becomes
a trivial modification of the proof of Theorem \ref{theosob}.  \qed

\section{An extension to the hard--core lattice gas (independent sets)}

\subsection{The hard core lattice gas}
A configuration $\eta\in\O:=\{0,1\}^{\Tree^b}$ is called an independent
set if no two adjacent vertices are occupied, \ie if $\eta_x\eta_y=0$
for every couple $x,y\in\Tree^b$ such that $d(x,y)=1$. We call $\wb\O$
the collection of all independent sets over the $b$--ary tree $\Tree^b$.
In the hard--core lattice gas model $\wb\O$ is the set of allowed
configurations and each such configuration $\h\in\wb\O$ is weighted with
the factor $\l^{|\h|}$ where $|\h|$ stands for the cardinality of $\h$,
\ie the number of occupied vertices in $\h$, and $\l> 0$ is the
so--called activity parameter.  To define the Gibbs measure we use local
specifications $\mu_A^\t$ obtained by setting
$$
\mu_A^\t(\eta) \propto \l^{|\eta_A|}\,,
$$
where $A$ is a finite subset of $\Tree^b$, $\tau,\eta\in\wb\O$ are
two allowed configurations such that $\tau_x=\eta_x$ for all $x\notin
A$, and $|\eta_A|=\sum_{x\in A}\eta_x$.  It is well known that the
hard--core lattice gas model undergoes a phase transition at the
critical activity $\l_c=b^b/((b-1)^{b+1})$ (see e.g.\ 
\cite{Spitzer,Kelly}).  For $\l\leq \l_c$ there is a unique phase
regardless of the boundary condition on the leaves, while for $\l>\l_c$
there are (at least) two distinct phases, corresponding to the
\emph{odd} and \emph{even} boundary conditions respectively.  The even
boundary condition $\tau^e$ is obtained by occupying all the vertices at
even depth from the root and letting all the rest unoccupied, \ie
$$
\tau^e_x=\begin{cases} 1 & d(x) \;\text{is even}\\
0  & d(x)\; \text{is odd}
\end{cases}
$$
The odd boundary condition $\tau^o$ is the complement
$\tau^o=1-\tau^e$.  We use the notation $\mu_{\ell}^e =
\mu_{T_\ell}^{\tau^e}$ for the Gibbs measure on the tree of depth $\ell$
with even boundary condition. Similarly $\mu_\ell^o$ denotes the Gibbs
measure with odd boundary conditions.  We also write
$\mu^e=\lim_{\ell\to\infty} \mu_\ell^e$ and $\mu^o=\lim_{\ell\to\infty}
\mu_\ell^o$.  When we need to emphasize the $\l$--dependence we shall
write $\mu_\l^e,\mu_\l^o$ in place of $\mu^e,\mu^o$.  Phase transition
is reflected by the fact that, when $\l>\l_c$, the probability of
occupation of the root differs for $\mu^e_\l$ and $\mu^o_\l$.

\subsection{The Glauber dynamics}
The hard--core Glauber dynamics is the Markov process with Markov
generator formally given by (\ref{generator}) with flip rates that are
reversible w.r.t.\ the hard--core lattice gas Gibbs measure.  As in the
Ising model we restrict for simplicity to the heat--bath dynamics given
by
\begin{equation}
c(\si^x)=\begin{cases} q_\l& \si_x=1\\
p_\l & \si^x\in\wb\O\,,\; \si_x=0\\
0 & \si^x\notin\wb\O\end{cases}
\;\quad\; q_\l:=\frac1{1+\l}\,,\; p_\l:=\frac\l{1+\l}\,.
\end{equation}
Here $\si\in\wb\O$ and $\si^x$ represents the configuration $\si$ with
the occupation number at $x$ inverted, \ie $(\si^x)_y=\si_y$, for all
$y\neq x$ and $(\si^x)_x=1-\si^x$.  In words, only transitions within
$\wb\O$ are allowed and the transition $\si\to\si^x$ occurs with rate
$p_\l=\l/(1+\l)$ if $x$ is vacant and with rate $q_\l=1/(1+\l)$ if $x$
is occupied. It is easily verified that detailed balance holds with this
choice of rates. Moreover, for any finite subset $A\subset \Tree^b$, for
any $\tau\in\wb\O$, the finite volume dynamics on $A$ with boundary
condition $\tau$ is ergodic and reversible w.r.t.\ the Gibbs measure
$\mu_A^\t$. As for the Ising Glauber dynamics we can use the spectral gap and
the logarithmic Sobolev constant to estimate the rate of convergence to
the stationary distribution $\mu_A^\t$. The corresponding definitions
are exactly the same as in (\ref{eq:cgapcsob}). An important result of
\cite{MaSiWe} is that the uniform bounds (\ref{eq:key-bound}) hold here
if we replace $\mu^+_\ell$ with $\mu^e_\ell$, \ie in the even phase one
has exponential decay to equilibrium for all values of $\l$. Of course,
the same holds for the odd phase.

\subsection{Attractivity}
It is essential for our approach that we can define a partial order on $\wb\O$
such that the hard--core lattice gas and its Glauber dynamics become \emph{attractive}.
Let us write $\Tree^b$ as the disjoint union of even and odd vertices, $\Teven$ and $\Todd$, where 
$\Teven:=\{x\in\Tree^b: \; d(x) \;{\rm is \;even}\}$ 
and $\Todd:=\{x\in\Tree^b: \; d(x) \;{\rm is \;odd}\}$. We define the following order on $\wb\O$:
\begin{equation}
\si\prec \h \;\iff\; \begin{cases}
\si_x\leq \h_x & x\in\Teven  \\
\si_x\geq \h_x & x\in\Todd
\end{cases}
\la{order}
\end{equation}
A function $f:\wb\O\to\bbR$ is called monotone increasing (decreasing)
if $\si\prec\h$ implies $f(\si)\leq f(\h)$ ($f(\si)\geq f(\h)$). We also
write $\mu \leq \nu$, for two measures on $\wb\O$, whenever $\mu(f)\leq
\nu(f)$ for every monotone increasing function $f$.  As in the Ising
model it is straightforward to construct an order--preserving global
path-wise coupling.  Let $\si_t^{\xi,A,\t}$ denote the hard--core Glauber
process at time $t$, with start in the configuration $\xi\in\wb\O$,
evolved in the region $A$ with boundary condition $\t\in\wb\O$.  We may
couple the processes $\bigl\{(\s_t^{\xi,A,\t})_{t\geq 0},\ A\sset
\Tree^b,\,\xi,\t\in\O\bigr\}$ such that the following relations hold:
for any $A\sset B\sset \Tree^b$, any $\xi\prec \xi'$, $\t\prec \t'$ and
$t\geq 0$
\begin{align}
  \s_t^{\xi,A,\t} &\prec \s_t^{\xi',A,\t'}  \\
\s_t^{\xi,A,\tau^o} &\prec \s_t^{\xi,B,\t} \prec \s_t^{\xi,A,\tau^e}  
\label{eq:mono200}
\end{align}
These relations also imply the following monotonicity properties of the Gibbs measures
and the associated FKG--property (see (\ref{mono1})--(\ref{eq:mono20})):
\begin{align}
&(i) \qquad \text{for any $A\sset \Tree^b$ the map $\h \mapsto
  \mu_A^\eta(f)$ is increasing;}
\label{mono1000}\\
&(ii) \qquad \text{$\mu_B^{\tau^e}\leq \mu_A^{\tau^e}$ whenever $A\sset B$.}
  \label{eq:mono2000}
\end{align}

\subsection{Results}
Replacing $\mu^+$ with $\mu^e$ we may define the sets $\O_\a$,
$\a\in(0,1)$, just as in Definition \ref{Omega alpha}.
It is not difficult to check that Lemma \ref{pro+} and therefore
Corollary \ref{incr} hold in the present setting as well as in the
Ising case. The same applies to Corollay \ref{weak} and Lemma \ref{mu+le}.

\smallskip

We need to introduce the hard--core analog of the Bernoulli measures
$\bbP_p$.  We call $\nu_{p,\l}$, $p\in(0,1), \l\geq 0$ the probability
measure on $\wb\O$ obtained as follows: we first assign occupation
numbers on $\Teven$ according to the Bernoulli(p) probability $\bbP_p$.
This gives a configuration $\h$ on $\Teven$. To obtain a legal
configuration (in $\wb\O$) we may now occupy only those vertices in
$\Todd$ that are at least at distance $3$ from $\h$. Call
$A_\h\subset\Todd$ this set of available vertices. Finally put
$\eta_x=1$ with probability $p_\l = \l/(1+\l)$ independently for every
$x\in A_\h$.

\smallskip

A simple coupling argument shows that, for every $\l>0$, $\nu_{p,\l}\geq
\mu_\l^e$ as soon as $p\geq p_\l$.  In particular, the argument of
Lemma \ref{mu+le} shows that $\nu_{p,\l}(\O_\a)=1$, for every $\a>0$,
for all $p\geq p_\l$. Our main result for the hard--core lattice gas is
stated as follows.

\begin{theorem}
\la{hard-corethm}
\medno
\begin{enumerate}[a)]
\item For every $b\geq 2$,
there exists $p<1$ such that
for all $\l\in(0,\infty)$  
we have $\nu(\O_\a)=1$, for some $\a=\a(\l,b)>0$, for
any initial distribution 
$\nu$ such that $\nu\geq
\nu_{p,\l}$\,.
\medno
\item For every $p>\frac12$, there exist $b_0\in\bbN$ and $\l_0\in(0,\infty)$ such that 
for $b\geq b_0$, $\l\geq \l_0$ we have $\nu(\O_\a)=1$, for some $\a=\a(\l,b,p)>0$,
for any initial distribution $\nu$ such that $\nu\geq\nu_{p,\l}$\,.
\end{enumerate}
\end{theorem}

\subsection{Sketch of proof of Theorem \ref{hard-corethm}}
Theorem \ref{hard-corethm} will be proved with the same arguments used
in the proof of Theorem \ref{theotree}. Below we point out the necessary (rather obvious)
modifications. 

The first observation is that in view of the monotonicity of $\O_\a$,
the domination $\nu_{p,\l}\geq \mu_\l^e$ for $p\geq p_\l$, 
Lemma \ref{mu+le} allows to 
replace the statements in the theorem by

\medno 
{\em
\begin{enumerate}[a*)]
\item{For every $b\geq 2$,
there exist $p<1$ and $\l_0<\infty$ such that
for all $\l\geq\l_0$  
we have $\nu_{p,\l}(\O_\a)=1$, for some $\a=\a(\l,b)>0$.}
\item{For every $p>\frac12$ there exist
$b_0$ and $\l_0$ such that 
for $b\geq b_0$, $\l\geq \l_0$
we have $\nu_{p,\l}(\O_\a)=1$
for some $\a=\a(\l,b,p)>0$.}
\end{enumerate}
} 

\medno

To repeat the argument of section \ref{main_proof} we need to introduce
the notion of the environment of obstacles. A realization of the
environment is described by $\o\in\wb\O$ with the following
interpretation: $x\in\Teven$ is called an \emph{obstacle} if $\o_x=0$ and
is said to be \emph{free} if $\o_x=1$.  Similarly $x\in\Todd$ is an
obstacle if $\o_x=1$ and is free if $\o_x=0$. Note that $x$ is an
obstacle in $\o$ iff $\o_x=\t^o_x$. As in the Ising case $\o$ determines
the tree $T(\o)$, \ie the largest connected component of free vertices
containing the root.  We write $\cB_\o$ for the set of $\t\in\wb\O$ such
that $\t_x=\t^o_x$ for every $x\notin T(\o)$. The hard--core model in a
given environment $\o$ is then obtained as before: for every
$A\subset\Tree^b$ we write $\mu_{A,\o}^\t$ for the measure $\mu_{A\cap
  T(\o)}^\t$, where $\t\in\cB_\o$.  The same reasoning applies to the
dynamics and we may use, as before, $\r_{t,\o}(\xi)$ for the expected
value at the root of the occupation variable under the dynamics 
$\si_{t,\o}^\xi$ among obstacles  
with starting configuration $\xi\in\wb\O$
%$\r_{t,\o}^{\ell,\t}(\xi)$ to denote
%the expected value at the root of the occupation variable under the dynamics
%in $T_\ell(\o)$ with starting configuration $\xi\in\wb\O$ and boundary
%condition $\t\in\cB_\o$. When $A=T(\o)$ we simply write
%and $\r_{t,\o}(\xi)$.

\smallskip

We then observe, as in 
(\ref{eq:mono4}) that  $$
\r_{t,\o}(\t^e)\leq \r_t(\o)\,,\qquad x\in\Teven\,,$$ 
where $\r_t(\o)$ denotes expectation at the root
w.r.t.\ the dynamics in infinite volume without obstacles
with starting configuration $\o$. 
We then define the environment $\o=\o(\h,\ell)$ 
as in (\ref{rule}), where of course the $+$ configuration is
replaced by $\t^e$. 
We now proceed exactly as in (\ref{main.1}).
Moreover, we may repeat the estimates of the three
terms there without modifications. 
What is crucial is that 
the technical estimates isolated in Claims 1 to 4 
can be established for the new setting.
A discussion of the point is given in the next subsection.

\subsection{Technical estimates}
To prove the Claims 1 to 4 for the hard--core model one needs to adapt
to the present setting the analysis developed in section \ref{recursive}. One defines 
the ratios $R$ and the associated weights $W$ in a similar way here, but 
the recursive relations involved in the proofs of the main estimates are model--specific
and require a separate investigation. We will not provide all the details
here since there is no truly new ingredient. However we give a sketch
of the
basic computations on the ratios $R$
to help the interested reader in reconstructing the needed claims. 

\smallskip\smallskip

\noindent
{\bf Estimates on $R$.}
Let $R$ be defined by
\begin{equation}
R(\o)=\frac{\mu^e_\o(\si_r=1)}{\mu^e_\o(\si_r=0)}
\la{rhc}
\end{equation}
A simple computation gives that if $\o_r=1$ we have %can calculate 
\begin{equation}
R(\o) = \l\,\prod_{i=1}^b\frac1{(1+R_{x_i}(\o))}
\la{iterhc}
\end{equation}
where $x_i$, $i=1,\dots,b$ denote the children of the root and $R_{x_i}(\o)$ is the 
corresponding ratio, given as usual by the rule $R_x(\o)=R(\theta_x\o)$ 
($\theta_x\o$ being
the environment shifted by $x\in\Tree^b$). 
The crucial estimate (\ref{R1}) is now replaced by 
\begin{equation}
\wt\nu_{p,\l}\left(R\leq \sqrt\l\right)\leq \d\,.
\la{Rhc}
\end{equation}
Here $\wt\nu_{p,\l}$ denotes the probability $\nu_{p,\l}$ conditioned to have 
$\o_r=1$ and $\d$ is a small parameter to be fixed at a later stage.
The following bound is the analogue of Lemma \ref{R_a} in the present setting.

\begin{lemma}
\la{R_ahc}
For any $\d>0$, $b\geq 2$, 
there exist $p_0<1$ and $\l_0<\infty$ such that
(\ref{Rhc}) holds for all $p\geq p_0$ and $\l\geq \l_0$.
\end{lemma}
\proof
We shall give the proof only in the case $b=2$.
For any integer $\ell$ we may define the ratios $R^\ell(\o)$ w.r.t.\ $\mu^e_{\ell,\o}$
as in (\ref{r2}). We set 
\begin{equation}
q_\ell = \wt\nu_{p,\l}\left(R^\ell \leq \sqrt\l\right)
\la{qellhc}
\end{equation}
Let $x_1,x_2$ denote the children of the root and call $y_1,y_2$ and $y_3,y_4$ 
the children of $x_1$ and $x_2$, respectively. Observe that the event $E$
that $\o_{x_i}=0$, $i=1,2$ and $\o_{y_i}=1$, $i=1,\dots,4$ has $\nu_{p,\l}$ 
probability at least $1-4(1-p)$ (since it suffices to occupy all $y_i$'s to
automatically free the $x_i$'s). Moreover for $\o\in E$ we have
\begin{equation}
R^\ell(\o)= \l\,\left(1+\frac{\l}{(1+R^\ell_{y_1}(\o))(1+R^\ell_{y_2}(\o))}\right)^{-1}
\left(1+\frac{\l}{(1+R^\ell_{y_3}(\o))(1+R^\ell_{y_4}(\o))}\right)^{-1}\,.
\la{qellhc2}
\end{equation}
Note that $R^\ell_{y_j}(\o))=R^{\ell-2}(\theta_{y_j}\o)$. 
Suppose that $R^\ell_{y_i}(\o) > \sqrt\l$, $i=1,2,3$. Then the above formula shows that
for $\l$ sufficiently large, the condition $R^\ell(\o)\leq \sqrt\l$ forces 
$R^\ell_{y_i}(\o)\leq 3$. Reasoning as in (\ref{qq2}) we see that 
\begin{equation}
q_\ell \leq 4(1-p) + 6 q_{\ell-2}^2 + 4 \wt\nu_{p,\l}\left(R^\ell_{y_1} \leq 3\right)\,.
\la{qellhc3}
\end{equation}
Using again (\ref{qellhc2}) we see that if there is only one $y_j$ 
with $R^\ell_{y_j}(\o)\leq \sqrt\l$ it is impossible to have $R^\ell(\o)\leq 3$.
It follows that 
\begin{equation} 
\wt\nu_{p,\l}\left(R^\ell_{y_1} \leq 3\right)
\leq 4(1-p) + 6 q_{\ell-4}^2 
\la{qellhc4}
\end{equation}
Putting these estimates together and using $q_{\ell-4}\leq q_{\ell-2}$ we see that
\begin{equation}
q_\ell \leq 12(1-p) + 30 q_{\ell-2}^2 \,.
\la{qellhc5}
\end{equation}
The conclusion now follows from (\ref{qellhc5}) just as in the case of (\ref{qq4}) 
because of the even boundary condition.
\qed

\section{Open problems}

We conclude by discussing an interesting open problem. Back to the Ising
case with $h=0$, let us take as initial distribution for the Glauber
dynamics the symmetric product measure $\bbP_{1/2}$ that for shortness
we denote by $\nu$. 

A first non trivial question is whether the law of the Glauber dynamics
$\nu P_t$ converges to a Gibbs measure as $t\to \infty$. In $\bbZ^d$ it
is well known that this is the case (see e.g. \cite{Li}) because $\nu
 P_t$ is translation invariant; unfortunately the Lyapunov function techniques
behind the proof do not seem to apply on the tree because of the large
boundary/volume ratio.

In the uniqueness region $\b\leq \b_0$ it is not difficult to check that
$\nu  P_t$ converges weakly to the unique Gibbs measure as $t\to \infty$.
More interesting is the interval $\b\in (\b_0,\b_1)$, where $\b_1$ is
the spin-glass transition point discussed in section
\ref{sec:Ising_on_trees}. Here the situation is more complicate due to
the presence of infinitely many extremal Gibbs states.

If we recall our first characterization of $\b_1$, it is not
unreasonable to conjecture that $\nu  P_t$ will
converge to the (extremal) free Gibbs measure $\mu^{\rm free}$. In fact,
if we
imagine that the single site Glauber dynamics is replaced by a block
heat bath dynamics as in \cite{BKMP} then, at least for small times,
each update of a block (say a large but finite subtree) replaces the
Bernoulli product measure $\nu$ inside the block with a finite Gibbs
measure close to $\mu^{\rm free}$. The case $\b\geq \b_1$ should be even more complex and one
can conceive that the dynamics and coarsening of clusters of spins with
opposite sign, present in the starting configuration, will play a
significant role as in the $\b=+\infty$ case \cite{H}. 

Although we have no clear answers to any of the above questions, we do
have some preliminary ``concentration of measures'' results that bring
some support to the conjectured behavior in the intermediate regime
$(\b_0,\b_1)$. 

First we show that for any
local function $f$ 
\begin{equation*}
\nu\Bigl(\h;\ | P_t(f)(\h)-\nu  P_t(f)|\geq \nep{-ct}\Bigr)\leq 
\nep{-c\nep{ct}}\,,%\quad c\ll 1
\end{equation*}
for some $c>0$.

In other words not too large fluctuations in the
starting configuration $\h$ are completely washed out by the
dynamics. In particular, for any local function $f$ which is \emph{odd}
w.r.t. a global spin flip, 
$$
\lim_{t\to \infty} P_t(f)(\h)=\mu^{\rm free}(f)=0\qquad \text{$\nu$--a.a. $\h$}\,.
$$
Secondly we derive a stability result that can be roughly formulated
as follows. Let $\tilde\nu$ be a perturbation of $\nu$ such that the relative
entropy between $\nu$ and $\tilde\nu$ restricted to the first $\ell$ levels
does not grow faster than $(b^{\ell})^\d$, $\d\ll 1$. Then for any local
function $f$ 
$$
\lim_{t\to \infty} |\tilde\nu  P_t(f) - \nu  P_t(f)|=0\,.
$$
We now formalize what we just said.

\begin{proposition}
\label{pro:free1}
For $h=0$ and $\b <\b_1$ there exists a positive constant $c>0$ such
that, for any function $f$ depending only on
finitely many spins and any $t\geq 0$ :
\begin{equation}
  \label{eq:free1}
\nu\Bigl(\h;\ | P_tf(\h)-\nu  P_tf|\ge
\nep{-ct}\Bigr)\leq \nep{-c_f\nep{ct}}\,.
\end{equation}
for a suitable constant $c_f >0$ depending on $f$.  
\end{proposition}
\proof
We are going to use standard Gaussian concentration bounds \cite{Ledoux}
for the measure
$\nu$ of the form:
\begin{equation}
  \label{eq:free01}
  \nu(\h;\ |F(\h)|\geq r) \leq \nep{-\frac{r^2}{2}}
\end{equation}
for any mean zero function $F$ with unitary Lipshitz norm 
$$
\|F\|^2_{\rm Lip}:= \sum_{x\in \Tree^b}\ninf{F(\h^x)-F(\h)}^2\,.
$$
Therefore (\ref{eq:free1}) follows if we can prove that for some $a=a(\b)>0$
\begin{equation}
  \label{eq:free001}
\| P_tf\|_{\rm Lip} \le
C_f\nep{- at}\,.
\end{equation}
for a suitable constant $C_f >0$ depending on $f$.
The basic tool is coupling
along the lines introduced in \cite{BKMP}. 
Recall that $\tanh(\b)< 1/\sqrt{b}$ for any $\b<\b_1$. Thus we can always choose $\l\in
\bigl(\tanh(\b),(b\tanh(\b))^{-1}\bigr)$ in  such a way that $b\l^2<1$. Given two
configurations $\h,\xi$ that differ in finitely many points, define
their weighted Hamming distance as
\begin{equation}
  \label{eq:free02}
  d_\l(\h,\xi)=\sum_x \l^{d(x)}\un_{\h_x\neq \xi_x}
\end{equation}
Then a key result of \cite{BKMP} combined with an
unpublished paper of Peres and Winkler (see section 4 of \cite{BKMP})
shows that under the natural coupling of the Glauber dynamics started at
$\h$ and $\xi$ 
\begin{equation}
  \label{eq:free03}
\bbE\bigl(d_\l(\s_t^\h,\s_t^\xi)\bigr)\leq C\nep{-ct}d_\l(\h,\xi)
\end{equation}
for suitable positive constants $C,\,c$.
Therefore
\begin{gather}
  \| P_tf(\h^x)- P_tf(\h)\|_{\infty}\leq \sum_y
   \ninf{f(\h^y)-f(\h)}\ \bbE\bigl(\un_{\s^{\h^x}_{t,y}\neq \s^{\h}_{t,y}}
   \bigr) \nonumber \\
\leq \Bigl(\sum_y \l^{-d(y)}\, \ninf{f(\h^y)-f(\h)} \Bigr)
   \nep{-ct}\l^{d(x)} = C_f \ \nep{-ct}\l^{d(x)}
\label{eq:free04}
  \end{gather}
Since $b\l^2 <1$, the sum over $x$ of the \emph{square} of the r.h.s.\
of (\ref{eq:free04}) converges and (\ref{eq:free001}) follows. \qed

\begin{corollary}
In the same setting as above, let $\tilde \nu$ be a probability measure on $\O$ and let $\nu_\ell,\ \tilde
\nu_\ell$ be the marginals on $\O_{T_\ell}$ of $\nu$ and $\tilde \nu$ respectively. Then there exists
$\d=\d(\b)$ such that, if $\ \Ent_\nu(\frac{d\tilde\nu_\ell}{d\nu_\ell})\le
b^{\d\ell}$ for all $\ell$, then     
\begin{equation}
  \label{eq:free2}
\lim_{t\to \infty} |\tilde\nu  P_t(f) - \nu  P_t(f)|=0\,.
\end{equation}
Moreover the limit is attained exponentially fast.
\end{corollary}
\proof Let $k=k(\b)$ be so large that, with $\ell=kt$, for any large
enough $t$
$$
\ninf{ P_t(f)-\nep{t\cL^{\rm free}_\ell}(f)}\leq
\nep{-t}\,,
$$ 
where $\cL^{\rm free}_\ell$ stands for the generator of the Glauber
dynamics in $T_\ell$ with free boundary conditions.
Standard results on finite speed of information propagation
show that such a $k$ exists (see e.g. \cite{Mar}). 
Let now $c_f,c$ be the
constants appearing in Proposition \ref{eq:free1} and let $A_t$ be the
set of configuration $\{\h;\ | P_t(f)(\h)-\nu  P_t(f)|\geq
\nep{-ct}\}$. 
Then, by setting $h_\ell:= \frac{d\tilde\nu_\ell}{d\nu_\ell}$,
\begin{gather}
  \label{eq:free3}
  |\tilde\nu  P_t(f) - \nu  P_t(f)| \leq 2\nep{-t}
  + |\tilde \nu_\ell \nep{t\cL^{\rm free}_\ell}(f)- \nu_\ell
  \nep{t\cL^{\rm free}_\ell}(f)|\nonumber\\
  = 2\nep{-t}
  +|\nu\bigl(\bigl[h_\ell -1\bigr] \nep{t\cL^{\rm
  free}_\ell}(f)\bigr)|
\nonumber\\
  \leq 4\nep{-t}+2\nep{-ct} + \nep{-c_f\nep{ct}} + \ninf{f}\nu\bigl(h_\ell\un_{A_t}\bigr)
\end{gather}
It remains to bound $\nu\bigl(h_\ell\un_{A_t}\bigr)$ and this is easily
accomplished using the entropy inequality together with Proposition
\ref{pro:free1} and our 
assumption on $\Ent_\nu(h_\ell)$.   
For any $\l>0$ 
\begin{gather}
  \label{eq:free4}
  \nu\bigl(h_\ell\un_{A_t}\bigr) \leq \frac 1\l
  \log\Bigl(\nu\bigl(\nep{\l\un_{A_t}} \bigr)\Bigr) + \frac 1\l
  \Ent_\nu(h_\ell)\\
\leq  \frac 1\l
  \log\Bigl(1+ (\nep{\l}-1)\nep{-c_f\nep{ct}}\Bigr) + \frac{b^{\d kt}}{\l}
\end{gather}
If we now choose $\l= \frac{1}{4}c_f\nep{ct}$ and $\d< \frac{c}{k\log b}$
we see that $\nu\bigl(h_\ell\un_{A_t}\bigr)$ tends to zero as $t\to
\infty$ exponentially fast. \qed

%%%%%%%%%%%%%%%%%%%%%%%%%%%%%%%%%%%%%%%%%%%%%%%%%%%%%%%%%%%%%%%%%%%
%%%%%%%%%%%%%%%%%%%%%%%%%%%%%%%%%%%%%%%%%%%%%%%%%%%%%%%%%%%%%%%%%%%
%%%%               REFERENCES                    %%%%%%%%%%%%%%%%%%
%%%%%%%%%%%%%%%%%%%%%%%%%%%%%%%%%%%%%%%%%%%%%%%%%%%%%%%%%%%%%%%%%%%
%%%%%%%%%%%%%%%%%%%%%%%%%%%%%%%%%%%%%%%%%%%%%%%%%%%%%%%%%%%%%%%%%%%

%%%%%%%%%%%%%%%%%%%%%%%%%%%%%%%%%%%%%%%%%%%%%%%%%%%%%%%%%%%%%%%%%%%
%%%%%%%%%%%%%%%%%         END OF PAPER       %%%%%%%%%%%%%%%%%%%%%%
%%%%%%%%%%%%%%%%%%%%%%%%%%%%%%%%%%%%%%%%%%%%%%%%%%%%%%%%%%%%%%%%%%%
\end{document}